\documentclass[12pt,a4paper]{article}
\usepackage[latin1]{inputenc}
\usepackage[english]{babel}
\usepackage[left=1.5cm,right=1.5cm,top=2cm,bottom=2.5cm]{geometry}
\author{Pieter Vanmechelen\thanks{KU Leuven, Department of Computer Science} \and Geert Lombaert\thanks{KU Leuven, Department of Civil Engineering} \and Giovanni Samaey\thanks{KU Leuven, Department of Computer Science. Corresponding author: \texttt{giovanni.samaey@kuleuven.be}}}
\date{}
\title{Comparison of random field discretizations for high-resolution Bayesian parameter identification in finite element elasticity}
\usepackage{amssymb}
\usepackage{mathtools}
\usepackage{amsthm}
\usepackage{graphicx}
\usepackage{listings}
\usepackage{enumitem}
\usepackage{multicol}
\usepackage{xcolor}
\usepackage{relsize}
\usepackage{subcaption}
\usepackage{tikz}
\usepackage{algorithm}
\usepackage{algorithmic}
\usepackage{url}
\graphicspath{{images/}}
\newtheorem{remark}{Remark}
\definecolor{KULdark}{RGB}{47,77,93}
\definecolor{KULred}{RGB}{147,53,53}
\definecolor{KULbright}{RGB}{29,141,176}
\definecolor{KULgreen}{RGB}{60,110,58}
\definecolor{KULmidlight}{RGB}{204,219,228}
\definecolor{KULverylight}{RGB}{231,238,242}
\definecolor{thesis}{RGB}{0,112,182}
\begin{document}

\maketitle

\begin{abstract}
We compare three random field discretization strategies for probabilistic identification of spatially varying material parameters in high-resolution finite element models. These strategies are (i) the Karhunen-Lo\`eve expansion, (ii) a wavelet expansion, and (iii) local average subdivision. The methods are assessed in the context of multilevel Markov chain Monte Carlo applied to plane stress elasticity with high-resolution displacement observations. Emphasis is placed on numerical efficiency, initialization cost, Markov chain mixing, and cost-to-error behaviour as the discretization resolution increases. While all approaches yield comparable posterior estimates, significant differences are observed in multilevel variance reduction and sampling efficiency. In particular, local average subdivision exhibits improved mixing and lower cost-to-error ratios at fine resolutions, despite its higher nominal parameter dimension. The results provide practical guidance for selecting stochastic field representations in uncertainty quantification in finite element simulations of heterogeneous materials.
\end{abstract}

\section{Introduction}
\label{ch:intro}

Bayesian inversion is a widely used approach to parameter estimation in many branches of science \cite{law_data_2015}. Assuming that we have access to a mathematical forward model of the system under study, that takes a parameter $E$ as input, we can use inversion strategies to estimate the value of $E$ conditional on observations of the model output. This is done in a Bayesian way: rather than simply finding a value for the parameter $E$ that optimally explains the observations $u$, a full posterior distribution $P(E|u)$ given the data $u$ is sought to enable uncertainty quantification. This is done using Bayes' rule in Equation \eqref{eq:Bayes}, which states that $P(E|u)$ is proportional to the product of (i) a prior distribution $\pi$ containing any information on $E$ known prior to the experiment and (ii) a likelihood function, which expresses how well a given realisation of $E$ matches with the observations $u$:

\begin{equation}
P(E|u) \simeq \mathcal{L}(u|E)\pi(E).
\label{eq:Bayes}
\end{equation}

In many applications, the parameter $E$ under study is a continuous function of space and the entire profile of $E$ is sought. One effective approach in this case is to model $E$ as a random field, commonly modelled as a Gaussian random field and subsequently transformed to achieve the desired characteristics of $E$. There is a vast literature on Gaussian random field modelling, with several known approaches \cite{liu_advances_2019}. For the context of this paper, three approaches will be compared and are thus of special interest. The first is the Karhunen-Lo\`eve (KL) expansion, which uses an eigenpair decomposition of the random field covariance function \cite{blondeel_p-refined_2020,ghanem_stochastic_1991,schwab_karhunen_2006}. The KL expansion has been extensively used in the context of multilevel inversion methods \cite{charrier_finite_2013,dodwell_hierarchical_2015,vanmechelen_multilevel_2025}. This approach is known to produce good posterior matching to the forward model output even for low truncation numbers in the eigenpair expansion \cite{uribe_bayesian_2020}.

The KL approach is based on global (eigen)functions of the random field covariance operator as a basis for the field. To contrast this, the two other methods under study are local expansions. The first alternative is the use of local average subdivision (LAS) \cite{fenton_simulation_1990,nuttall_parallel_2011}. This method uses a set of increasingly local averages to maximise local information, but suffers more strongly from the curse of dimensionality as its dimension is inherently tied to the resolution of the FE grid. The second alternative approach is a wavelet-based expansion \cite{bachmayr_representations_2016}, which aims to compromise between the two extremes.

To compare the these three approaches in a Bayesian setting, we will make use of the Markov chain Monte Carlo (MCMC) methodology. Here, we construct the posterior distribution from a Markov chain of samples that are proposed according to some proposal distribution $q$, and subsequently accepted or rejected, based on their relative posterior scores. These samples can then be used to construct Monte Carlo (MC) estimates for relevant quantities of interest.

The random field approach to modelling $E$ typically results in very high-dimensional parameter representations. We additionally consider cases where the forward model is discretised on a high-resolution grid. These two aspects combined result in a challenging MCMC problem. We will use the multilevel MCMC algorithm \cite{dodwell_hierarchical_2015,vanmechelen_multilevel_2025} with preconditioned Crank-Nicolson (pCN) proposals \cite{cotter_mcmc_2013}, which is especially suited for this high-dimensional, high-resolution case. This is but one strategy, and many alternatives exist such as delayed acceptance MCMC \cite{christen_markov_2005, lykkegaard_multilevel_2020}, sequential MC \cite{hoang_analysis_2020,latz_multilevel_2018}, transitional MCMC \cite{angelikopoulos_x-tmcmc_2015,ching_transitional_2007}, or dimension-independent likelihood-informed MCMC \cite{cui_multilevel_2024,cui_likelihood-informed_2014}.

The goal of this paper is to provide a computational comparison between the three approaches in a specific Bayesian inversion setting. We focus on a problem that arises in high-resolution FE simulations with uncertain material fields, as commonly encountered in computational mechanics. The forward model consists of solving the linear elasticity equations in the plane stress case on a rectangular domain. The model is discretised on finite element (FE) grids with varying resolutions in the multilevel hierarchy. The inferred parameter is the spatial variation of the Young's modulus $E$, conditioned on observations of a static displacement field $u$ under a prescribed load.

The article is structured as follows. In Section \ref{ch:model} the linear elasticity model is explained. Section \ref{ch:representations} shows methods of constructing the Young's modulus field for all three approaches. Section \ref{ch:mcmc} explains the MCMC approach, starting from a Metropolis-Hastings implementation and then the multilevel MCMC method. In Section \ref{ch:numerics}, a comparison of MCMC results is performed for the three representations in a numerical setup in a high-resolution, high-dimensional setting including high-resolution data. Finally, Section \ref{ch:conclusion} contains the main conclusions of this study.

All code used to generate the results in this article is open-source, and can be found at \url{https://gitlab.kuleuven.be/numa/public/mlmcmc-estimations.git}.

\section{Model outline}
\label{ch:model}

We will first discuss the underlying model used within this article, which is the equations of linear elasticity. Subsection \ref{sec:linear_elasticity} details the nature of these equations, while Subsection \ref{sec:random_field} focuses on the idea of modelling the parameter under study as a random field.

\subsection{Linear elasticity model}
\label{sec:linear_elasticity}

\subsubsection{Equations of linear elasticity}

The model studied in this paper consists of the 2D equations of linear elasticity in the plane stress case. This is a set of equations that links the spatially varying 2D Young's modulus $E(x,y)$ and Poisson's ratio $\tilde{\nu}$ of the underlying material to the 2D displacement field $\pmb{u}(x,y)$ of the system when it is subjected to a set of body forces $\pmb{F}_{\text{body}}(x,y)$. By introducing the Lam\'e parameters $\{\lambda,\mu\}$, it is possible to combine the linear elasticity equations into one elliptic PDE of the following form:

\begin{subequations}
\begin{align}
&-\nabla\cdot\left[\lambda(\nabla\cdot \pmb{u})I + \mu(\nabla \pmb{u} + (\nabla \pmb{u})^T)\right] = \pmb{F}_{\text{body}} 
\\
&\lambda = \dfrac{E\tilde{\nu}}{(1+\tilde{\nu})(1-2\tilde{\nu})}, \qquad \mu = \dfrac{E}{2(1+\tilde{\nu})}.
\end{align}
\label{eq:linear_elasticity}
\end{subequations}

\begin{figure}
\centering
\resizebox{\textwidth}{!}{
\begin{tikzpicture}
\fill[color={rgb,255:red,231; green,238; blue,242}] (0,0) rectangle (16, 4);
\foreach \x in {1,...,31}
	{\draw (\x/2,0) -- (\x/2,4);
	\fill[color={rgb,255:red,60; green,110; blue,58}] (\x/2, 0) circle (0.07);
	\fill[color={rgb,255:red,60; green,110; blue,58}] (\x/2, 4) circle (0.07);}
\foreach \y in {0,...,8}
	{\draw (0, \y/2) -- (16, \y/2);}
\draw[color={rgb,255:red,60; green,110; blue,58}, line width=0.5mm] (0,0) rectangle (16, 4);
\draw[ultra thick] (0,-0.6) -- (0, 4.6);
\draw[ultra thick] (16,-0.6) -- (16, 4.6);
\foreach \y in {0,...,21} 
	{\draw (-1, \y/5+0.4) -- (0, \y/5-0.6);
	\draw (16, \y/5+0.4) -- (17, \y/5-0.6);}
\foreach \x in {0,...,3} 
	{\draw (-1, 0.2-\x/5) -- (-0.2-\x/5, -0.6);
	\draw (-0.8+\x/5, 4.6) -- (0, 3.8+\x/5);
	\draw (16, -0.4+\x/5) -- (16.2+\x/5, -0.6);
	\draw (16.2+\x/5, 4.6) -- (17, 3.8+\x/5);}
\foreach \x in {13,...,26}
	{\draw[->, ultra thick] (0.375*\x+0.1875, 4.5) -- (0.375*\x+0.1875, 4);}
\end{tikzpicture}%
}
\caption{Clamped beam setup, with a load placed on the center third of the top edge. Observations are made at all nodes along the top and bottom edges of the beam (green lines). This sketch shows a resolution of $32\times 8$ finite elements.}
\label{fig:beam}
\end{figure}
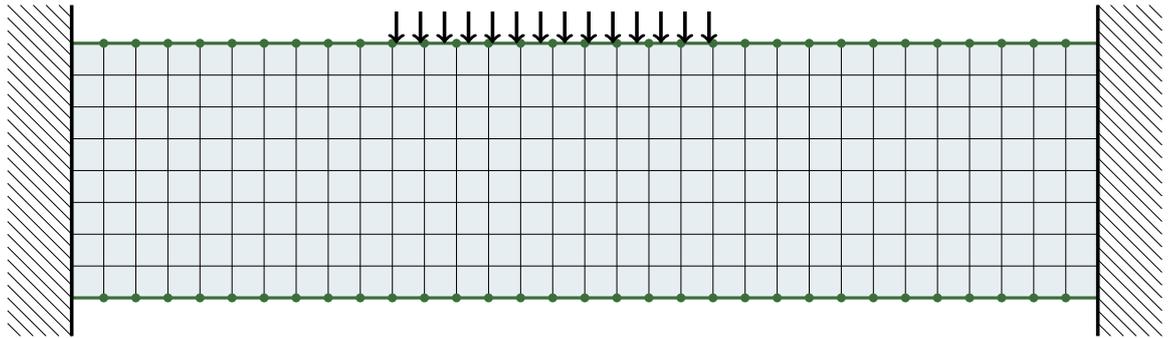

In this equation, $I$ is the identity matrix. For simplicity, in the context of this paper we assume the Poisson's ratio $\tilde{\nu}$ to be constant and known throughout the system, such that for a given $\pmb{F}_{\text{body}}$ the only input stiffness parameter to be modelled in this equation is the Young's modulus $E$. Within this paper, we will apply the linear elasticity equations to the case of a 2D rectangular domain, representing a beam as sketched in Figure \ref{fig:beam}. The beam is clamped at the left and right edges, giving Dirichlet boundary conditions $\pmb{u}=0$ along both edges. Additionally, an external downward line force is applied to the top edge of the domain. The horizontal position, length and magnitude of this force will be detailed in Section \ref{ch:numerics}. For practical calculations, we solve the PDE in Equations \eqref{eq:linear_elasticity} on a finite element (FE) grid. Linear square Lagrangian finite elements are chosen to discretize the model parameters. We will denote the grid-discretized version of $E$ as $E_{\ell}$ depending on the level $\ell$ of resolution.

\subsubsection{Synthetic data}

The output of the linear elasticity model is the displacement field $\pmb{u}$. We assume to have observations of this field at the top and bottom edges of the domain. While this output quantity is a continuous field, in practice observations are discretized, which are considered here to be of high resolution. This high resolution is modelled in the FE approach by assuming that observations are made at each node on the top and bottom edges of the grid, which is thus constructed at a high resolution itself. In practice, such high-resolution measurements can be obtained using techniques such as Digital Image Correlation \cite{helfrick_3d_2011}. 

In this article, we will use synthetic observations, generated by solving Equations \eqref{eq:linear_elasticity} for a known ground truth profile of $E$. To simulate measurement error, these observations are perturbed by noise, which is assumed to be independent and normally distributed. Hence the observations used in the model under study are collected in the vector $\pmb{u}^{\text{obs}}$, which is calculated as

\begin{equation}
\pmb{u}^{\text{obs}} = \mathcal{P}(\pmb{u}) + \sigma_F^2\mathcal{N}(0,I),
\label{eq:noise}
\end{equation}

where $\mathcal{P}$ is the projection operator onto the nodes on the top and bottom edges of the FE grid. In total there are $N_{\text{obs}}$ such nodes, yielding $\pmb{u}^{\text{obs}} \in \mathbb{R}^{2N_{\text{obs}}}$.

\subsection{The Young's modulus as a random field}
\label{sec:random_field}

We are interested in modelling the spatially varying Young's modulus $E$ in Equations \eqref{eq:linear_elasticity} outlined above. We first construct the Young's modulus field as a Gaussian random field $g(x,\omega)$ and subsequently transform it to fit the desired properties of the Young's modulus. A transformation of $g$ that is regularly used in the context of linear elasticity results in a Gamma random field as follows \cite{blondeel_p-refined_2020}

\begin{equation}
E(x,\omega) = \mu\gamma^{-1}\left[\kappa, \dfrac{\Gamma(\kappa)}{2}\left(1+\text{erf}\left(\dfrac{g(x,\omega)}{\sqrt{2}}\right)\right)\right].
\label{eq:transformation}
\end{equation}

Here $\Gamma$ denotes the Gamma function, $\gamma^{-1}$ the inverse of the lower incomplete gamma function and erf the error function. $\mu$ and $\kappa$ are the scale and shape parameters of the Gamma field \cite{blondeel_p-refined_2020}.

With the transformation defined, modelling the Young's modulus field reduces to modelling a Gaussian random field. In what follows, the computational domain of the field is $D = [0,4]\times[0,1]$, corresponding to the beam dimensions of Figure \ref{fig:beam}. The random field additionally requires a covariance kernel, detailing the covariance between two points $x,y\in D$. Here, we choose a Mat\'ern covariance kernel on the Gaussian field $g$, which is well-studied in the context of random field modelling and a popular choice for many applications \cite{lord_introduction_2014}. Denoting its variance as $\sigma^2$, length scale as $\lambda$ and smoothness as $\nu$, the Mat\'ern kernel takes the form

\begin{equation}
C(x,y) = \dfrac{\sigma^2}{2^{\nu-1}\Gamma(\nu)}\left(\sqrt{2\nu}\dfrac{\|x-y\|_p}{\lambda}\right)^{\nu}K_{\nu}\left(\sqrt{2\nu}\dfrac{\|x-y\|_p}{\lambda}\right),
\label{eq:matern}
\end{equation}

where $K_{\nu}$ is the modified Bessel function of the second kind. Here we will consider isotropic covariance functions for which $p=2$.

\section{Three different random field representations}
\label{ch:representations}

From a computational mechanics perspective, the choice of random field representation determines not only approximation accuracy but also solver cost, scalability, and algorithmic robustness in stochastic finite element simulations. Section \ref{ch:model} showed how to model the Young's modulus as a transformed normal field $g$. The following subsections detail three possible approaches to sampling $g$: using the Karhunen-Lo\`eve basis in Subsection \ref{sec:kl}, using wavelets in Subsection \ref{sec:wavelet} and using local average subdivision (LAS) in Subsection \ref{sec:las}.

\subsection{Karhunen-Lo\`eve approach}
\label{sec:kl}

\begin{figure}[h]
\centering
\hspace*{-2cm}
\includegraphics[width=1.2\textwidth]{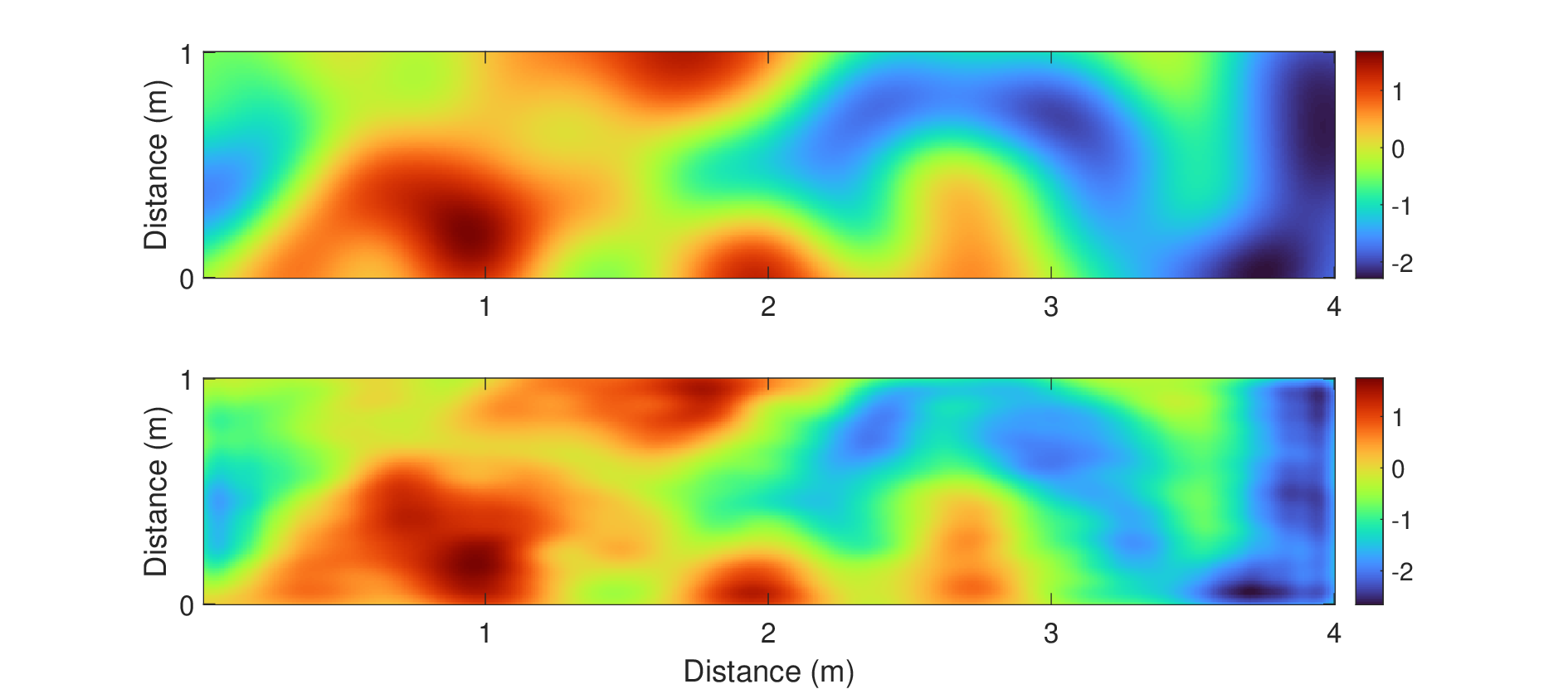}
\caption{Gaussian random field generated using the KL expansion, for different values of truncation number $m$. Top: $m=64$. Bottom: $m=1024$.}
\label{fig:kl_convergence}
\end{figure}

The Karhunen-Lo\`eve (KL) approach starts with finding pairs of eigenvalues $\lambda_m$ and eigenfunctions $b_m(x)$ of the covariance operator of $g$, obtained by solving the eigenvalue problem

\begin{equation}
\int_D C(x,y)b_m(y)dy = \lambda_mb_m(x).
\label{eq:eigenvalue_problem}
\end{equation}

With these eigenpairs obtained, we can decompose the Gaussian field $g$. This provides a useful series expansion for $g$, where the terms are ordered by decreasing eigenvalue $\lambda_m$:

\begin{equation}
g(x,\omega) = \mathbb{E}[g(\cdot,\omega)] + \sum_{m=1}^{\infty}\sqrt{\lambda_m}\xi_m(\omega)b_m(x),
\label{eq:kl}
\end{equation}

where $\xi_m \sim \mathcal{N}(0,1)$ for all $m$. This expansion can be truncated at an arbitrary numer of terms. As the eigenpairs can be precomputed, we only need to generate the realisations of $\xi_m$ in a sampling algorithm. As such, the KL method can easily be used in very high-dimensional cases, by carrying the expansion above to a high truncation value of $m$. The KL basis has been well-studied in inverse problems governed by forward PDEs \cite{cliffe_multilevel_2011,lord_introduction_2014}. It has also been used specifically in linear elasticity forward problems, for example by \cite{blondeel_p-refined_2020}. The main appeal of the KL expansion is that it minimises the mean squared error of a truncated expansion of a random field using orthonormal basis functions \cite{lord_introduction_2014}.

Figure \ref{fig:kl_convergence} shows one realisation of a Gaussian field with Mat\'ern covariance kernel with $\sigma = 1$, $\nu = 1.5$, $\lambda = 0.5$ generated using the KL expansion, for two different values of the truncation number $m$.

\subsection{Wavelet approach}
\label{sec:wavelet}

\begin{figure}[h]
\centering
\hspace*{-2cm}
\includegraphics[width=1.2\textwidth]{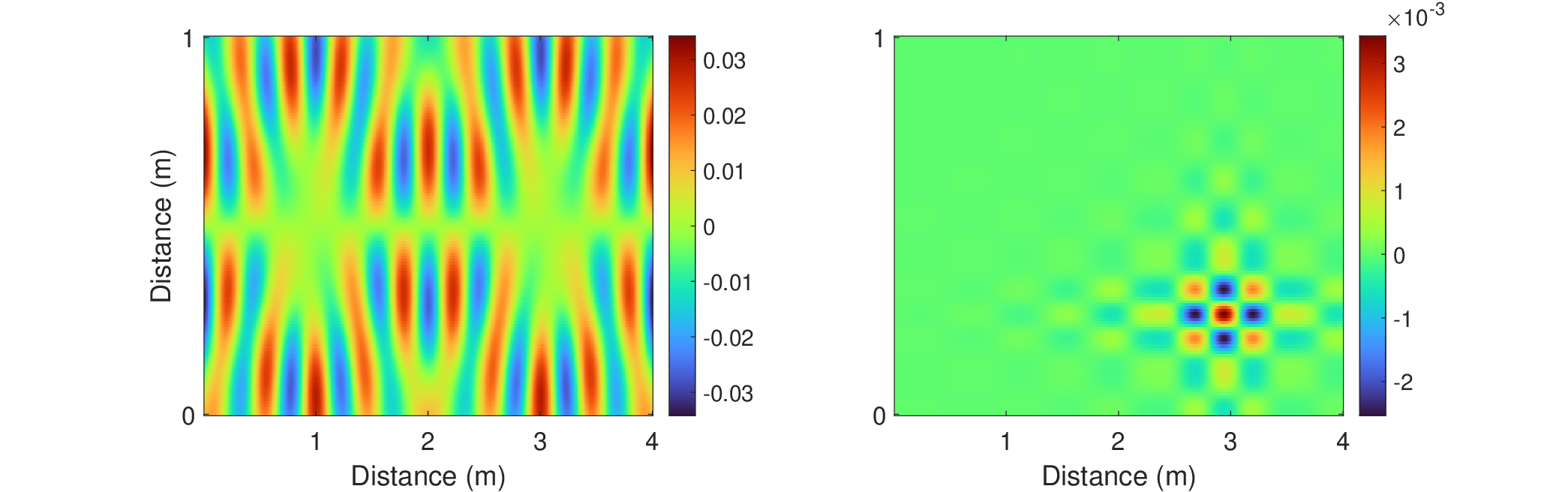}
\caption{One basis function each from KL (left) and wavelet (right) expansions. Note that the KL basis function is $L^2$-normalised, whereas the wavelet is not as its norm decay is directly encoded in the construction of the basis functions.}
\label{fig:kl_wavelet}
\end{figure}

An alternative approach to the global basis functions of the KL expansion is to use local basis functions. To this extent, the next expansion under study is a wavelet-based approach \cite{bachmayr_representations_2016}. To illustrate the difference, Figure \ref{fig:kl_wavelet} shows one basis function from both approaches. While the KL basis function is oscillating throughout the domain, the wavelet oscillation is confined to part of the domain.

The remainder of this subsection outlines how to construct the wavelet basis for a Gaussian random field with covariance function $C$ on a general $d$-dimensional domain $D$. For a more in-depth discussion, we refer the reader to \cite{bachmayr_representations_2016}.
To generate a wavelet basis for a Gaussian random field on $D$, the first step is to embed $D$ in a $d$-dimensional torus $\mathbb{T} = [-\gamma,\gamma]^d \supset D$. The periodic nature of $\mathbb{T}$ allows for the straightforward construction of a wavelet basis $\{\tilde{b}_m\}_{m\in\mathbb{N}}$ on $L^2(\mathbb{T})$ \cite{bachmayr_representations_2016}. In the context of this paper, we choose the basis of Meyer wavelets. Next, we construct the Fourier expansion of each wavelet basis function $\tilde{b}_m$:

\begin{equation}
\tilde{b}_m = (2\gamma)^{-d/2}\sum_{n\in\mathbb{Z}^d} c_n(\tilde{b}_m)e_n,
\end{equation}

where

\begin{equation}
e_n(z) = (2\gamma)^{-d/2}e^{i\pi n\cdot z/\gamma}\quad \text{and} \quad c_n(\cdot) = \int_{\mathbb{T}} \cdot e^{-i\pi n\cdot z/\gamma}dz, \quad n \in \mathbb{Z}^d.
\end{equation}

The third step is to construct a covariance function on $\mathbb{T}$ that agrees exactly with $C$ over $D$. This is done by defining a sufficiently smooth and even cutoff function $\phi: \mathbb{R}^d \to \mathbb{R}$, for which $\phi\vert_D = 1$ and $\phi(x)=0$ for $\|x\|_{\infty} > 2\gamma - \text{diam}(D)$. Let $\text{diam}(D) = \delta$, then we choose

\begin{equation}
\phi(x) = \dfrac{\zeta\left(\frac{\kappa-|x|}{\kappa-\delta}\right)}{\zeta\left(\frac{\kappa-|x|}{\kappa-\delta}\right) + \zeta\left(\frac{|x|-\delta}{\kappa-\delta}\right)} \qquad \text{and } \zeta(x) = \left\lbrace \begin{array}{ll}
\exp(-x^{-1})& \text{if }x > 0 \\ 0 & \text{if }x \leq 0,
\end{array} \right.
\end{equation}

with $\kappa = 2\gamma - \delta$. The covariance function is then defined as

\begin{equation}
C^{\mathbb{T}} = \sum_{n\in\mathbb{Z}^d}C(z+2\gamma n)\phi(z+2\gamma n).
\end{equation}

In this notation, $C$ is viewed as a function of only one variable, i.e. $\|x-y\|_p$ in Equation \eqref{eq:matern}. A Gaussian field on $\mathbb{T}$ can then be expanded as

\begin{equation}
g^{\mathbb{T}}(x,\omega) = \sum_{m=1}^{\infty} \xi_m(\omega)b^{\mathbb{T}}_m(x),
\label{eq:wavelet_expansion}
\end{equation}

where again $\xi_m \sim \mathcal{N}(0,1)$ for all $m$, and \cite{bachmayr_representations_2016}

\begin{equation}
b^{\mathbb{T}}_m(x) = (2\gamma)^{-d/2}\sum_{n\in\mathbb{Z}^d}\sqrt{c_n(C^{\mathbb{T}})}c_n(\tilde{b}_m)e_n.
\label{eq:waveletsum}
\end{equation}

\begin{figure}[h]
\centering
\hspace*{-2cm}
\includegraphics[width=1.2\textwidth]{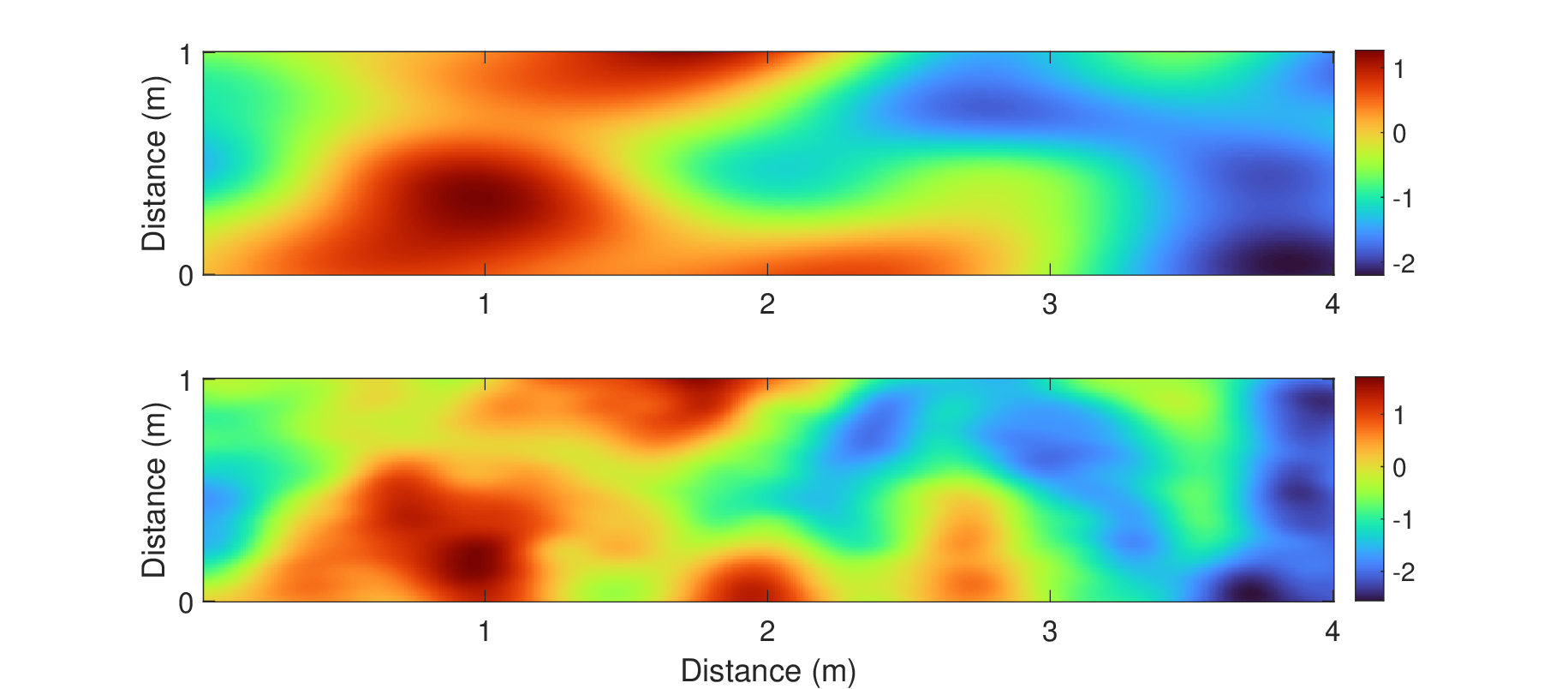}
\caption{Gaussian random field generated using the wavelet expansion, for different values of truncation number $m$. Top: $m=64$. Bottom: $m=1024$.}
\label{fig:wavelet_convergence}
\end{figure}

The desired Gaussian field $g$ on $D$ is simply the restriction of $g^{\mathbb{T}}$ to $D$. Notice that, unlike Equation \eqref{eq:kl}, Equation \eqref{eq:wavelet_expansion} does not explicitly contain the expectation of $g$. In the wavelet basis, this value is encoded in the wavelet scaling function, which is represented by $b_1^{\mathbb{T}}(x)$. Figure \ref{fig:wavelet_convergence} shows the same Gaussian random field as before, using the wavelet basis and truncated at 64 and 1024 modes, illustrating that this method yields very similar results to the KL method.

\subsection{Local average subdivision approach}
\label{sec:las}

The third approach of constructing a Gaussian field is local average subdivision (LAS) \cite{fenton_simulation_1990,nuttall_parallel_2011}. Unlike the previous two approaches, this method does not aim to construct a set of basis functions to use in an expansion like Equations \eqref{eq:kl} and \eqref{eq:wavelet_expansion}. Instead, it constructs the field from a sequence of grid refinements, where at each iteration the value inside each cell of the domain is the local average of the random field within the region defined by this cell.

In the domain $D = [0,4]\times [0,1]$ discretised with square finite elements outlined so far, this means starting from 4 square elements at the initial iteration $k=0$, which are directly sampled from a Gaussian distribution with covariance based on the pairwise distances between the elements. At each subsequent iteration $k$, every cell is then divided into 4 smaller squares. The value within these new cells is determined by the value in the parent cell, as well as its direct (including diagonal) neighbours. The updating procedure can be expressed as a linear system, denoting by $g^{k}_{i,j}$ the mean value of $g$ in the cell at iteration $k$ and position $(i,j)$:

\begin{subequations}
\begin{align}
g^{k+1}_{2i,2j} &= c_{11}\xi^{k+1}_{i,j,1} + \sum_{m=-1}^1\sum_{n=-1}^1 a^{k}_{m,n,1}g^{k}_{i+m,j+n}
\\
g^{k+1}_{2i,2j-1} &= c_{21}\xi^{k+1}_{i,j,1} + c_{22}\xi^{k+1}_{i,j,2} + \sum_{m=-1}^1\sum_{n=-1}^1 a^{k}_{m,n,2}g^{k}_{i+m,j+n}
\\
g^{k+1}_{2i-1,2j} &= c_{31}\xi^{k+1}_{i,j,1} + c_{32}\xi^{k+1}_{i,j,2} + c_{33}\xi^{k+1}_{i,j,3} + \sum_{m=-1}^1\sum_{n=-1}^1 a^{k}_{m,n,3}g^{k}_{i+m,j+n}
\\
g^{k+1}_{2i-1,2j-1} &= 4g^{k}_{i,j} - g^{k+1}_{2i,2j} - g^{k+1}_{2i,2j-1} - g^{k+1}_{2i-1,2j}.
\end{align}
\label{eq:LAS_update}
\end{subequations}

The coefficients $a$ and $c$ depend only on the covariance function $C$ (and resolution level $k$), and as such can be precomputed before any sampling routine is performed, see \cite{nuttall_parallel_2011} for details on their construction. Hence we only need to calculate the coefficients $\xi$ for a given realisation of the random field and as before, each coefficient is standard normally distributed. The fourth equation in this set is defined to ensure that the local average value of the parent cell is maintained after refinement. For a detailed implementation strategy of the 2D LAS method, we refer the reader to \cite[Chapter 3]{nuttall_parallel_2011}.

\begin{figure}[h]
\centering
\hspace*{-2cm}
\includegraphics[width=1.2\textwidth]{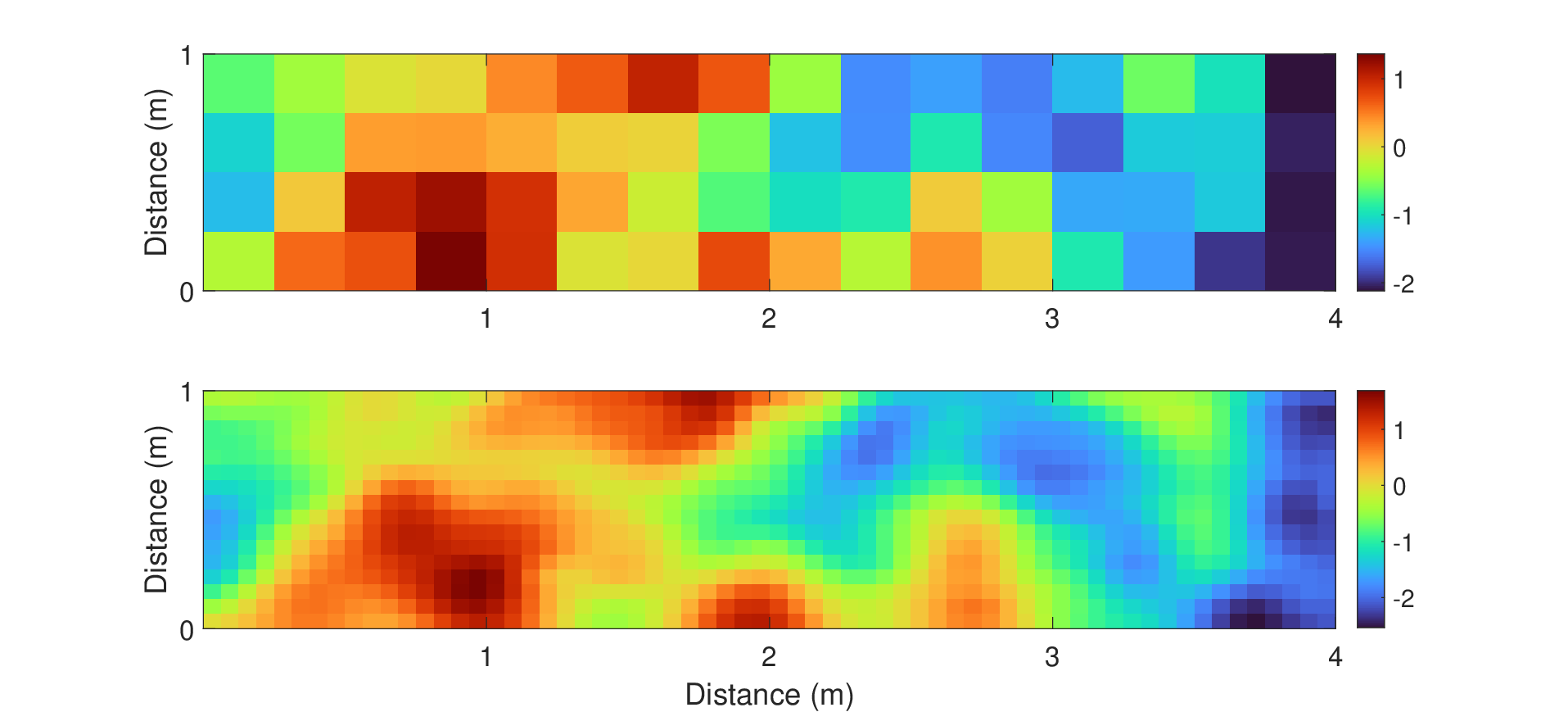}
\caption{Gaussian random field generated using the LAS method, for different number of refinement iterations. Top: max iteration $\ell=3$ (64 terms). Bottom: max iteration $\ell=5$ (1024 terms).}
\label{fig:las_convergence}
\end{figure}

The LAS method requires as many coefficients $\xi$ as cells used in the grid to sample a random field. As such, it samples the field using a significantly higher parameter dimension than either KL or wavelet expansions, as shown in Figure \ref{fig:las_convergence}. This does not necessarily lead to a loss of efficiency, as we will illustrate in Section \ref{ch:numerics}. Additionally, it offers some situational advantages over the two other approaches: its use of local averages means that it is very well suited for adaptive mesh refinement in multilevel strategies. Additionally, especially for very high resolution grids, the initial cost of the LAS is considerably lower than that of the two other methods, see Figure \ref{fig:costs}.

\begin{remark}
In its current form, the LAS method relies on a rectangular geometry of the domain and finite elements \cite{liu_advances_2019}. Extensions to this framework typically consist of embedding a grid into a rectangular domain or warping a rectangular grid to obtain more complex shapes \cite{griffiths_probabilistic_2004}. It would be of value to extend the division analysis of the LAS approach to triangular elements to enable more general geometries. However, such a treatment is beyond the scope of this thesis but should be relatively straightforward.
\end{remark}

\subsection{Convergence rates}

We now provide bounds on the $\mathcal{L}^2$-convergence of the three different representations to the true underlying field. Let $E_m^{\text{KL}}$ denote the field obtained by projecting a given realisation of $E$ onto the KL basis and truncating the expansion after $m$ terms. It has been shown that, for a $d$-dimensional Mat\'ern covariance with smoothness parameter $\nu > 1$ and under mild conditions for the transformation of Equation \eqref{eq:transformation}, the following inequality holds up to a constant independent of $m$: \cite{vanmechelen_multilevel_2025}

\begin{equation}
\|E_m^{\text{KL}}-E\|_{\mathcal{L}^2(D)} \lesssim m^{-\nu/d}.
\label{eq:kl_convrate}
\end{equation}

Using the wavelet approach, a similar bound has been proven for Gaussian fields $g$ \cite{bachmayr_representations_2016}. Following the reasoning in \cite{vanmechelen_multilevel_2025}, the same bound can also be obtained for the transformed field $E_m^{\text{Wavelet}}$ using the transformation of Equation \eqref{eq:transformation}, where $E_m^{\text{Wavelet}}$ denotes the gaussian field projected onto the wavelet basis, truncated after $m$ terms and fed through the transformation of Equation \eqref{eq:transformation}:

\begin{equation}
\|E_m^{\text{Wavelet}}-E\|_{\mathcal{L}^2(D)} \lesssim m^{-\nu/d}.
\end{equation}

For the LAS approach, the truncated quantity $E_m^{\text{LAS}}$ for $m=2^{dk}$ with positive integer $k$ is the approximation to $E$ after $k$ iterations in the LAS refinement. For other values of $m$, we define $E_m^{\text{LAS}}$ as the result of carrying out the updating strategy of Equation \eqref{eq:LAS_update} from left to right, top to bottom. This labelling is fairly arbitrary, but poses no issue as in practice only values for $m = 2^{dk}$ are used. As the LAS method consists of approximating a field by its local averages, the convergence rate is the classical rate for constant finite elements \cite{brenner_mathematical_2010}:

\begin{equation}
\|E_m^{\text{LAS}}-E\|_{\mathcal{L}^2(D)} \lesssim m^{-\frac{1}{d}}.
\label{eq:las_convrate}
\end{equation}

\begin{figure}[h]
\centering
\includegraphics[width=0.7\textwidth]{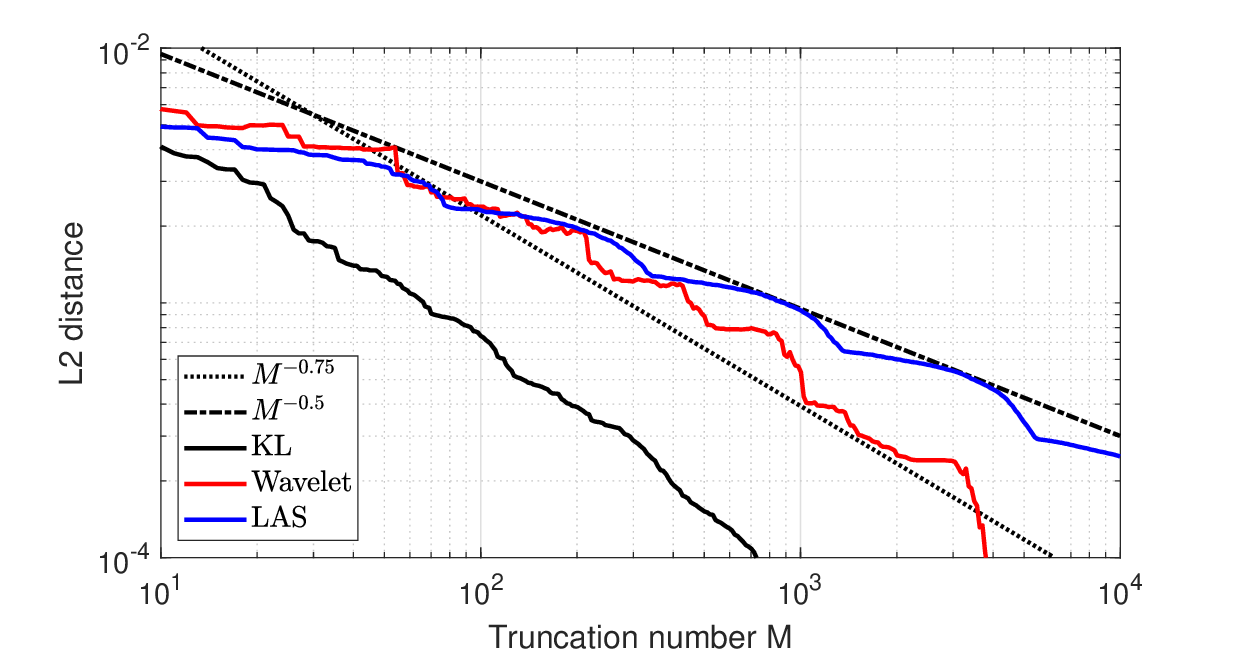}
\caption{$\mathcal{L}^2$-norms of the distance between truncated expansions of the field shown in Figures \ref{fig:kl_convergence}-\ref{fig:las_convergence} and the field itself, for each of the three representations discussed above.}
\label{fig:truncations}
\end{figure}

Figure \ref{fig:truncations} shows the convergence for each of the methods. We first generate the field shown earlier in this section using the wavelet approach for $m=5000$. We then project the field onto the KL basis functions and calculate its local means for the LAS terms. Each of the three methods is then truncated after $m$ terms and the $\mathcal{L}^2$ distance to the true field is calculated. As this field has a Mat\'ern covariance with smoothness $\nu = 1.5$, the KL and wavelet expansions converge considerably faster than the LAS approximation.

\section{Markov Chain Monte Carlo methods}
\label{ch:mcmc}

We will solve the inverse problem of recovering Young's modulus field profiles from displacement observations in a Bayesian framework. The current gold standard of this Bayesian approach is the Markov Chain Monte Carlo (MCMC) methodology \cite{law_data_2015}. In this section, we illustrate two MCMC methods: the basic Metropolis-Hastings algorithm \cite{hastings_monte-carlo_1970}, a simple and widely used MCMC method and the multilevel Markov Chain Monte Carlo method \cite{dodwell_multilevel_2019,vanmechelen_multilevel_2025}, which adapts Metropolis-Hastings for more efficient use in the case of high-resolution FE models with high-resolution observations.

\subsection{Bayesian inversion}
\label{sec:bayesian}

The goal of most MCMC methods is to sample from a posterior distribution $P(E|\pmb{u}^{\text{obs}})$ of the quantity of interest $E$ conditional on observations $\pmb{u}^{\text{obs}}$. As before $E$ is the value of the Young's modulus throughout the domain and $\pmb{u}^{\text{obs}}$ is the vector of displacements at the observation points. The main underlying idea of MCMC is Bayes' rule

\begin{equation}
P(E|\pmb{u}^{\text{obs}}) = \dfrac{\mathcal{L}(\pmb{u}^{\text{obs}}|E)\pi(E)}{P(\pmb{u}^{\text{obs}})},
\end{equation}

which states that the posterior distribution is proportional to the product of a prior distribution $\pi$ on $E$ and a likelihood function $\mathcal{L}$, which compares output of the forward model for a given value of $E$ to the available observations $\pmb{u}^{\text{obs}}$. The denominator in this expression is called the evidence, and is typically both infeasible and unnecessary to calculate \cite{law_data_2015}.

The use of Bayes' rule allows MCMC methods to sample from the posterior distribution by only needing to evaluate the likelihood and prior. By making the assumption of normally distributed measurement noise as in Equation \eqref{eq:noise}, the likelihood function takes the form

\begin{equation}
\mathcal{L}(\pmb{u}^{\text{obs}}|E) = \exp\left(-\dfrac{\|\mathcal{F}(E) - \pmb{u}^{\text{obs}}\|^2}{2\sigma_F^2}\right),
\label{eq:like}
\end{equation}

with the forward operator $\mathcal{F}$ defined as the concatenation of the solution operator of the PDE and the projection of this solution onto the measurement locations. As all three methods of sampling $E$ discussed in the previous section sample standard normal coefficients, the prior $\pi$ in this case is just a standard normal prior on the coefficient level.

\subsection{Metropolis--Hastings}
\label{sec:mh}

The Metropolis-Hastings algorithm is an MCMC algorithm that works by sequentially proposing a new sample $E'$ from some distribution $q$ conditional on the previous sample $E^n$ in the Markov chain and either accepting or rejecting it based on some criterion. The choice of proposal distribution is therefore an important factor in this algorithm, and several common choices exist, such as a random walk, where 

\begin{equation}
q^{\text{RW}}(\pmb{\xi}'|\pmb{\xi}^n) = \pmb{\xi}^n + \beta\mathcal{N}(0,I),
\label{eq:randomwalk}
\end{equation}

for some $\beta \geq 0$. The parameter $\pmb{\xi}^n$ is the set of coefficients in any of the sampling methods that represent the sample $E^n$. As these are equivalent, we will often write $q(E'|E^n)$ instead. A different choice is the preconditioned Crank-Nicholson (pCN) proposal, which is well-suited for high-dimensional problems such as random field modelling \cite{cotter_mcmc_2013}. This distribution takes the form

\begin{equation}
q^{\text{pCN}}(\pmb{\xi}'|\pmb{\xi}^n) = \sqrt{1-\beta^2} \pmb{\xi}^n + \beta\psi,
\label{eq:pcn}
\end{equation}

where $\psi$ is drawn from the prior distribution of $\pmb{\xi}$ in the case of a zero-mean Gaussian prior as for the methods outlined above, and $\beta \in [0,1]$. The pCN method works well in the case of Gaussian priors, with the main drawback that its use for general non-Gaussian priors is less straightforward. The Metropolis-Hastings algorithm for a general choice of proposal distribution $q$ is outlined in Algorithm \ref{alg:mh}.

\begin{algorithm}
\caption{Metropolis-Hastings}
\label{alg:mh}
\begin{algorithmic}
\STATE{Choose initial state $E^0 \in \mathbb{R}^{M}$}
\FOR{$n=0\ldots N-1$}
\STATE{Generate proposal move $E' \sim q(\cdot|E^n)$}
\STATE{Evaluate acceptance rate
\[
\alpha = \min\left\lbrace1,\dfrac{P(E'|\pmb{u}^{\text{obs}})q(E^n|E')}{P(E^n|\pmb{u}^{\text{obs}})q(E'|E^n)}\right\rbrace
\]}
\STATE{With probability $\alpha$ accept proposal and set $E^{n+1}=E'$}
\STATE{Otherwise reject proposal and set $E^{n+1}=E^n$}
\ENDFOR
\RETURN{Markov chain of elements $\{E^n\}$}
\end{algorithmic}
\end{algorithm}

\subsection{Multilevel Markov Chain Monte Carlo}
\label{sec:mlmcmc}

The Metropolis-Hastings algorithm requires, for each proposed sample $E'$ in the chain, solving the forward model $\mathcal{F}(E')$ to calculate the acceptance rate $\alpha$ of the sample. As such, the algorithm is ill-suited for cases where the forward model is expensive, because fewer samples can be generated, leading to poor posterior estimates. The multilevel Markov Chain Monte Carlo algorithm \cite{dodwell_multilevel_2019} tackles this problem by exploiting a hierarchy of cheaper surrogate models. This hierarchy is defined on levels $\ell=0,\ldots,L$ where the model itself is evaluated on the finest level $L$. For ease of presentation, we will add a subscript $\ell$ to all level-dependent variables and operations, to emphasize their dependency on the level.

In addition to increasing model fidelity between levels, at finer levels the model output is also typically of higher dimension. Using expansion-based methods such as the KL and wavelet expansions, this higher dimension simply manifests as a higher truncation number for the expansion. In the LAS approach, the dimension of the parameter is inherently tied to the grid resolution, so the parameter dimension automatically increases as the grid is refined. 

MCMC methods create posterior estimates for quantities of interest $Q$ by means of Monte Carlo estimates. In the multilevel case, this has been extended to multilevel Monte Carlo estimates \cite{giles_multilevel_2008,giles_multilevel_2015} by use of a telescopic sum

\begin{equation}
\widehat{Q}^{MLMC} = \dfrac{1}{N_0}\sum_{n=B_0+1}^{B_0+N_0}Q_0(E_0^n) + \sum_{\ell=1}^L\dfrac{1}{N_L}\sum_{n=B_{\ell}+1}^{B_{\ell}+N_{\ell}}\left[Q_{\ell}(E_{\ell}^n)-Q_{\ell-1}(\mathcal{E}_{\ell-1}^n)\right],
\label{eq:telescopic_sum}
\end{equation}

where $\{\mathcal{E}_{\ell-1}^n\}_{n=0}^{N_{\ell}}$ is a subsampled version of the chain at level $\ell-1$. The reason for using two chains will be detailed later in this section. The $B_{\ell}$ in this sum are called the burn-in lengths at each level. These are incorporated to discard the first few samples of the chain, as MCMC methods only asymptotically sample from the posterior distribution, and introducing a burn-in period decreases bias in the MCMC estimator \cite{law_data_2015}.

The telescopic sum in equation \eqref{eq:telescopic_sum} is significantly cheaper to evaluate than a direct Monte Carlo estimate of $Q_L$ based purely on samples at the finest level $L$. This is because each individual term in this sum is very cheap to evaluate: for low levels $\ell$, a coarse model is used so many samples $E_{\ell}^n$ can be generated. For higher levels and under some minor assumptions, the correction terms $Q_{\ell}-Q_{\ell-1}$ decrease as $\ell$ increases. As such, the variance of the estimator in the $\ell$-th term of this sum gets smaller, and fewer samples are needed to reach a given target threshold of the Monte Carlo error for this estimator. This decrease in samples needed typically more than compensates for the increased cost per sample at higher levels.

\begin{algorithm}[h]
\caption{Multilevel Markov chain Monte Carlo}
\label{alg:mlmcmc}
\begin{algorithmic}
\STATE{Choose initial state $E_L^0 \in \mathbb{R}^{M_L}$. For every $\ell \in 0,\ldots, L$ the initial state $E_{\ell}^0$ on level $\ell$ is given by the first $M_{\ell}$ modes of $E_L^0$. Choose subsampling rates $\tau_{\ell}$ for $\ell=1,\ldots,L$. Set level iteration numbers $n_{\ell}=0$ for $\ell=0,\ldots,L$}
\FOR{$\ell = 0,\ldots, L$}
\FOR{$n_{\ell} = 0,\ldots,N_{\ell}$}
\IF{$\ell = 0$}
\STATE{Generate new sample $E_0^{n_0+1}$ using Algorithm 1}
\ELSE
\STATE{Generate proposal coarse modes by subsampling chain on level $\ell - 1$: $(E_{\ell}^C)' = E_{\ell-1}^{n_{\ell}\tau_{\ell}}$}
\STATE{Generate proposal fine modes using pCN proposal distribution: $(E_{\ell}^F)' \sim q_{\ell}(\cdot|E_{\ell}^{n_{\ell},F})$}
\STATE{Set $E_{\ell}' = [(E_{\ell}^C)', (E_{\ell}^F)']$ and calculate
\[
\alpha_{\ell} = \min\left\lbrace1, \dfrac{\mathcal{L}_{\ell}(\pmb{u}^{\text{obs}}|E_{\ell}')\mathcal{L}_{\ell-1}(\pmb{u}^{\text{obs}}|E_{\ell}^{n,C})}{\mathcal{L}_{\ell}(\pmb{u}^{\text{obs}}|E_{\ell}^n)\mathcal{L}_{\ell-1}(\pmb{u}^{\text{obs}}|(E_{\ell}^C)')}\right\rbrace
\]}
\STATE{With probability $\alpha_{\ell}$ accept proposal and set $E_{\ell}^{n_{\ell}+1} = E_{\ell}'$,}
\STATE{Otherwise reject and set $E_{\ell}^{n_{\ell}+1} = E_{\ell}^{n_{\ell}}$}
\STATE{$n_{\ell} \leftarrow n_{\ell}+1$}
\ENDIF
\ENDFOR
\ENDFOR
\RETURN{Markov chains of elements $\{E^n_{\ell}\}$ at different discretisation levels}
\end{algorithmic}
\end{algorithm}

The multilevel MCMC algorithm \cite{dodwell_multilevel_2019} is concerned with generating the level-dependent chains $\{E_{\ell}\}_{n=0}^{N_{\ell}}$. In practice, these chains can be generated in parallel, but for simplicity the following outline shows a version of the algorithm where chains are generated sequentially. At the coarsest level $\ell=0$, the Metropolis-Hastings algorithm is used to generate the chain. The subsequent levels in the multilevel algorithm are generated in a different fashion. Because the dimension $M_{\ell}$ of the parameter increases between levels, different proposal methods are used to generate a proposal for the first $M_{\ell-1}$ modes than for the next $M_{\ell}-M_{\ell-1}$ modes, which are called coarse and fine modes, respectively. 

The coarse modes are generated by subsampling the previous level chain: by defining a subsampling rate $\tau_{\ell}\in\mathbb{N}$, a subsampled chain $\{\mathcal{E}_{\ell-1}^{n}\}_{n=0}^{N_{\ell}} = \{E_{\ell-1}^{n\tau_{\ell}}\}_{n=0}^{N_{\ell}}$ is constructed. The proposal for the coarse modes is then simply the next element in the chain $\mathcal{E}_{\ell-1}$. Care should be taken when choosing subsampling rates $\tau_{\ell}$: if $\tau_{\ell}$ is too small, the proposed samples are too highly correlated, leading to biased estimators. On the other hand, if $\tau_{\ell}$ is too large, the method loses computational efficiency as many more coarse samples must be calculated to generate a single proposal for the fine-level chain \cite{dodwell_multilevel_2019}. 

A proposal for the fine modes is generated by means of a proposal move starting from the previously accepted fine modes. To combat the curse of dimensionality, samples for the fine modes are typically generated using the pCN proposal distribution. Under this distribution, the acceptance rate takes the form

\begin{equation}
\alpha_{\ell} = \min\left\lbrace1, \dfrac{\mathcal{L}_{\ell}(\pmb{u}^{\text{obs}}|E_{\ell}')\mathcal{L}_{\ell-1}(\pmb{u}^{\text{obs}}|E_{\ell}^{n,C})}{\mathcal{L}_{\ell}(\pmb{u}^{\text{obs}}|E_{\ell}^n)\mathcal{L}_{\ell-1}(\pmb{u}^{\text{obs}}|(E_{\ell}^C)')}\right\rbrace.
\end{equation}

The acceptance rate in this equation is simply a ratio of likelihood scores for the proposed and previously accepted samples at the current and previous coarser levels, where $E_{\ell}^C$ denotes the coarse modes of $E_{\ell}$. Unlike equation \eqref{eq:like}, this likelihood is now level-dependent and takes the form

\begin{equation}
\mathcal{L}_{\ell}(\pmb{u}^{\text{obs}}|E_{\ell}) = \exp\left(-\dfrac{\|\mathcal{F}_{\ell}(E_{\ell})-W_{\ell}(\pmb{u}^{\text{obs}})\|^2}{2\mathcal{N}_{\ell}\sigma_F^2}\right).
\label{eq:levellike}
\end{equation}

Here $\mathcal{N}_{\ell}$ is the amount of nodes on the observed edges (top and bottom) of the $\ell$-th level mesh. The level-dependency of this likelihood is manifested in several places. First, the Young's modulus field $E$ is represented using a level-dependent number of modes $M_{\ell}$. Second, this field $E_{\ell}$ is fed through a forward model with a level-dependent mesh resolution. Third, to account for the discrepancy between mesh resolution of the coarse model output $\pmb{F}_{\ell}(E_{\ell})$ and the observations $\pmb{u}^{\text{obs}}$ which have a resolution corresponding to the finest level $L$, a weighting function $W_{\ell}:\mathbb{R}^{\mathcal{N}_L}\to\mathbb{R}^{\mathcal{N}_{\ell}}$ is applied to the observations to obtain a level-dependent surrogate set of observations that can be compared to the output of $\mathcal{F}_{\ell}$. Finally, a factor $\mathcal{N}_{\ell}$ is added to the denominator to account for the varying dimensionality of the terms in the numerator. This form of the likelihood function ensures consistency across all levels in the multilevel algorithm \cite{vanmechelen_multilevel_2025}.

The multilevel MCMC algorithm as outlined above is shown in Algorithm \ref{alg:mlmcmc}.

\section{Numerical comparison of parameter representations}
\label{ch:numerics}

The KL, wavelet and LAS methods will now be compared in a simple model problem of Bayesian inversion for a plane stress problem of linear elasticity, for two different ground truth cases. In Subsection \ref{sec:experiment_setup}, the details of the experiment setups are discussed. Subsection \ref{sec:posterior} contains a comparison of the posterior estimates obtained with each method for both experiments. In Subsection \ref{sec:cost_comparison}, the methods are compared with respect to their computational efficiency.

\subsection{Experiment setup}
\label{sec:experiment_setup}

The inverse problem under study is concerned with recovering the spatially varying value of the Young's modulus $E$ in a rectangular domain, conditional on the displacement observations made as outlined in Section \ref{ch:model}. The computational domain is again taken as $D = [0,4]\times[0,1]$, representing a 4 m $\times$ 1 m beam. 

It is beneficial to have a one-dimensional quantity of interest $Q$ summarising a given realisation of the spatial profile of $E$, to quantify the convergence rates of the MCMC algorithm. Here we choose the mean value of the integrated bending stiffness along the beam axis:

\begin{equation}
Q = \dfrac{1}{D_x} \int_0^{D_x}dx \int_0^{D_y} E(x,y) (y-y_0)^2 dy,
\label{eq:moment}
\end{equation}

where $x$ and $y$ are the horizontal and vertical coordinates of the domain, $D_x$ and $D_y$ are the length and height of the beam and the neutral line $y_0$ is the line of barycenters of $E$ in vertical slices of the domain:

\begin{equation}
y_0(x) = \dfrac{\int_0^{D_y} E(x,y)y\ dy}{\int_0^{D_y} E(x,y)\ dy}.
\end{equation}

In the context of this paper, synthetic data will be used. It is generated by assuming a certain ground truth profile of the Young's modulus and solving Equations \eqref{eq:linear_elasticity} with this input. We consider two different profiles here, detailed below. The other material parameters used to solve Equations \eqref{eq:linear_elasticity} are based on values for reinforced concrete: the Poisson's ratio $\tilde{\nu}$ is assumed to be 0.25, below the limit of incompressible regimes occurring at $\tilde{\nu} = 0.5$. The material density is assumed to be $2500$ kg/m$^3$. The external force in Equations \eqref{eq:linear_elasticity} is chosen as a 1 kN line load applied to the middle third of the top edge of the beam. 

The discretised domain on which the multilevel MCMC estimates are made is a grid consisting of 128 $\times$ 32 square finite elements. To avoid a so-called ``inverse crime", the synthetic data for both ground truths is generated at a higher resolution of 256 $\times$ 32 elements. The grid hierarchy used for all simulations in the multilevel algorithm consists of 4 levels $\ell = 0,\ldots,3$ with $2^{\ell+4}\times 2^{\ell+2}$ elements each.

\begin{figure}[h]
\centering
\includegraphics[width=\textwidth]{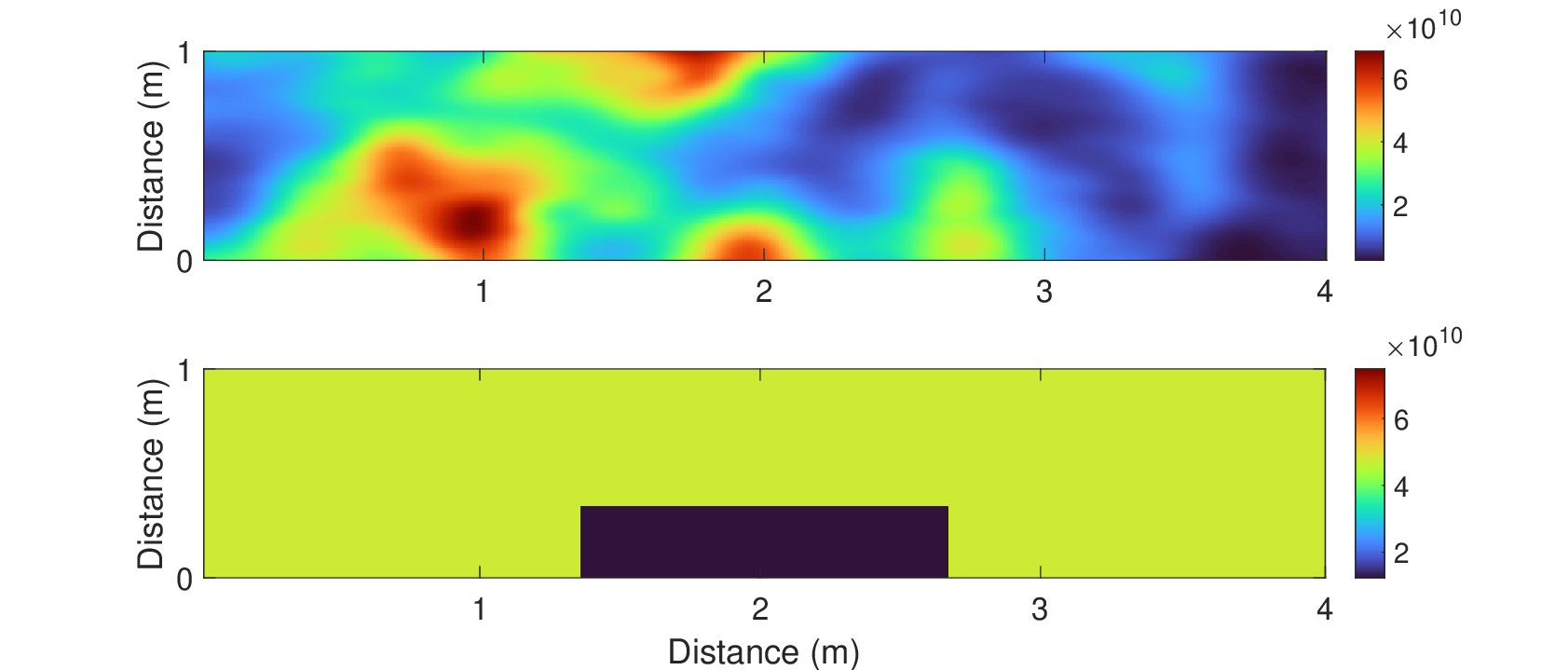}
\caption{Top: Ground truth profile used for the first experiment. Generated by sampling 5000 independent standard normal variables and using them in the wavelet expansion. Bottom: Ground truth profile used for the second experiment. Values are 12 GPa in the dark blue region, 47 GPa elsewhere.}
\label{fig:groundtruth}
\end{figure}

For the first experiment, we consider a ground truth for $E$ obtained by sampling from the truncated wavelet expansion using 5000 coefficients which are sampled from a standard normal prior. The first 1024 coefficients of these were used to generate the profile shown in Figure \ref{fig:wavelet_convergence}. This profile is subsequently transformed to a gamma random field according to Equation \eqref{eq:transformation}. The resulting field is shown in the top plot of Figure \ref{fig:groundtruth}. The bottom plot contains the ground truth profile for the second experiment, where we consider a piecewise constant ground truth profile. Here $E$ is assumed to take on a value of 47 GPa everywhere except for the bottom middle third of the beam, where a damaged region is assumed to be. Within this region $E$ is reduced to a value of 12 GPa.

The data obtained by solving Equations \eqref{eq:linear_elasticity} with these profiles as input are subsequently perturbed by noise according to Equation \eqref{eq:noise}. Due to the relatively low 4:1 aspect ratio of the domain, the resulting displacements are very small, on the order of $10^{-6}$ m. To obtain a controlled and representative signal-to-noise ratio, which is the relevant quantity for Bayesian inference as it determines the width of the likelihood function, we therefore add noise on the order of $10^{-9}$ to $10^{-8}$ m. While such absolute noise levels are smaller than those associated with realistic displacement measurement techniques such as Digital Image Correlation \cite{helfrick_3d_2011}, the resulting signal-to-noise ratio is comparable to that encountered in practical applications. The specific magnitudes of the signal and noise are thus chosen for visualization purposes, as the low aspect ratio facilitates clearer representation of the interior Young's modulus field, whereas the qualitative and quantitative behaviour of the proposed inference strategies depends on the signal-to-noise ratio rather than on the absolute displacement scale.

For the prior-sampled ground truth we perform the experiment with a value of $\sigma_F = 10^{-8}$, which corresponds to a measurement error on the order of $1\%$. This value enables a reasonably large step size in the MCMC routines and thus provides good mixing to compare the efficiency of the three approaches. The piecewise constant ground truth experiment uses a smaller value of $\sigma_F = 10^{-9}$ to perform a comparison between the methods in the case of a more strongly peaked likelihood.

For each of the three representations, the multilevel MCMC algorithm is used with pCN proposal moves. Here we choose subsampling rates $\tau_{\ell} = 100$ for $\ell = 1$ and $\tau_{\ell} = 5$ for $\ell > 1$. The pCN proposal moves use a step size of $\beta = 10^{-1}$ for the prior-sampled experiment and $\beta = 10^{-2}$ for the piecewise constant experiment, to account for the more peaked likelihood in this case. These values were found to provide (approximately) minimal value of the autocorrelation of the coarsest chain in the multilevel hierarchy for all three representations. For simplicity the step sizes at finer levels $\ell>0$ are chosen equal to the one at the coarsest level. We carry out the multilevel MCMC method for each representation on 10 parallel chains, each generating $1.5 \cdot 10^7$ coarsest level samples.

For ease of implementation, each method is thus compared using the same number of samples. To validate that this leads to similarly converged chains and that sufficiently many samples are taken, for each representation the Gelman-Rubin statistic \cite{gelman_inference_1992} is calculated of the chain estimating the quantity of interest in Equation \eqref{eq:moment} at all levels. The estimates for both experiments are shown in Table \ref{tab:gelmanrubin}.

\begin{table}
\centering
\begin{tabular}{| c | l | l | l | l |} 
\multicolumn{5}{c}{\textbf{Prior-sampled ground truth}} \\[2mm]
\hline
 & $\ell=0$ & $\ell=1$ & $\ell=2$ & $\ell=3$ \\ 
\hline
KL & 1.0002 & 1.0008 & 1.0006 & 1.0007 \\ 
\hline
Wavelet & 1.0002 & 1.0004 & 1.0006 & 1.002 \\ 
\hline
LAS & 1.00003 & 1.00005 & 1.00007 & 1.00002 \\ 
\hline
\multicolumn{5}{c}{} \\[2mm]
\multicolumn{5}{c}{\textbf{Piecewise constant ground truth}} \\[2mm]
\hline
 & $\ell=0$ & $\ell=1$ & $\ell=2$ & $\ell=3$ \\ 
\hline
KL & 1.001 & 1.003 & 1.004 & 1.004 \\ 
\hline
Wavelet & 1.006 & 1.01 & 1.006 & 1.01 \\ 
\hline
LAS & 1.001 & 1.003 & 1.004 & 1.005 \\ 
\hline
\end{tabular}
\caption{Estimates of the Gelman-Rubin statistic for the quantity of interest at each level in the multilevel hierarchy, for all three representations and both experiment setups.}
\label{tab:gelmanrubin}
\end{table}

The Gelman-Rubin statistic estimates the ratio between the within-chain and between-chain variances for multiple chain setups. Values close to 1 point to well-converged chains (typically $<1.1$ or $<1.01$ is aimed for). For all cases at all levels, we conclude based on this statistic that the chains have converged well, as all values are below 1.01, often considerably so.

The multilevel MCMC method not only requires a hierarchy of grid resolutions, but also a hierarchy of parameter dimensions. To determine appropriate dimensions for the three sampling methods we follow the traditional assumption on the bound of the mean squared error $e(\widehat{Q}_{L\{N_{\ell}\}}^{\text{MC}})^2$ for a quantity of interest $Q$ in the multilevel algorithm:

\begin{equation}
e(\widehat{Q}_{L\{N_{\ell}\}}^{\text{MC}})^2 \lesssim (L+1)\sum_{\ell=0}^L \dfrac{\mathbb{V}(Y_{\ell})\cdot \text{IAT}_{\ell}}{N_{\ell}} + \left(\mathbb{E}(Q_L) - \mathbb{E}(Q)\right)^2.
\label{eq:mse}
\end{equation}

The first terms in this assumption contain the sample variance $\mathbb{V}(Y_{\ell})$ of the correction terms, the integrated autocorrelation times $\text{IAT}_{\ell}$ of the chains and the sample sizes $N_{\ell}$. This sum forms the sampling error which will be discussed in more detail later. The second term is the bias error from approximating the true quantity of interest on a finite resolution grid in a finite-dimensional parameter space. We aim to make this bias error approximately equal in each of the three representations. As mentioned, we consider the same grid resolution in each case. The truncation $m_{\text{LAS}}$ of the LAS parameter naturally follows from this resolution and is taken to be $2^{2\ell+6}$ at each level. To find appropriate truncation numbers $m_{\text{KL}}$ for the KL and $m_{\text{W}}$ for the wavelet expansion, we look to Figure \ref{fig:truncations} and choose values such that

\begin{equation}
\|g_{m_{\text{KL}}}^{\text{KL}} - g\|_{\mathcal{L}^2(D)} \approx \|g_{m_{\text{W}}}^{\text{Wavelet}} - g\|_{\mathcal{L}^2(D)} \approx \|g_{m_{\text{LAS}}}^{\text{LAS}} - g\|_{\mathcal{L}^2(D)}.
\end{equation}

Table \ref{tab:trunc} lists the resulting truncations. It can be remarked that instead of choosing truncations based on the convergence of $g$, one should instead look at the convergence of $Q$. However, the convergence of $Q$ in terms of truncation is considerably less monotonous than that of $g$, making this approach less straightforward. Additionally, $Q$ follows the same asymptotic convergence rate as $g$ \cite{dodwell_hierarchical_2015} which implies $m_{\text{KL}},m_{\text{W}} \simeq m_{\text{LAS}}^{1/\nu}$. The rates shown in Table \ref{tab:trunc} roughly follow this rule for $\ell\geq 1$.

\begin{table}[h]
\centering
\begin{tabular}{c|ccc}
Level & $m_{\text{KL}}$ & $m_{\text{W}}$ & $m_{\text{LAS}}$ \\
\hline
0 & 16 & 64 & 64 \\
1 & 40 & 256 & 256 \\
2 & 100 & 500 & 1024 \\
3 & 250 & 1400 & 4096
\end{tabular}
\caption{Truncation numbers for each of the expansions at the four levels of the multilevel algorithm.}
\label{tab:trunc}
\end{table}

Finally, the samples used in both experiments model Gamma random fields with Mat\'ern covariance obtained using the transformation of Equation \eqref{eq:transformation}. The Mat\'ern parameters used are $\lambda=0.5$, $\nu=1.5$ and $\sigma=1$, which are the same parameters used to generate the prior-sampled ground truth case.

\subsection{Posterior estimates}
\label{sec:posterior}

We first compare the qualitative use of the three methods for making posterior predictions. We first show the results for the prior-sampled experiment and then the results for the piecewise constant experiment.

\subsubsection{Prior-sampled ground truth}

\begin{figure}[h]
\centering
\includegraphics[width=\textwidth]{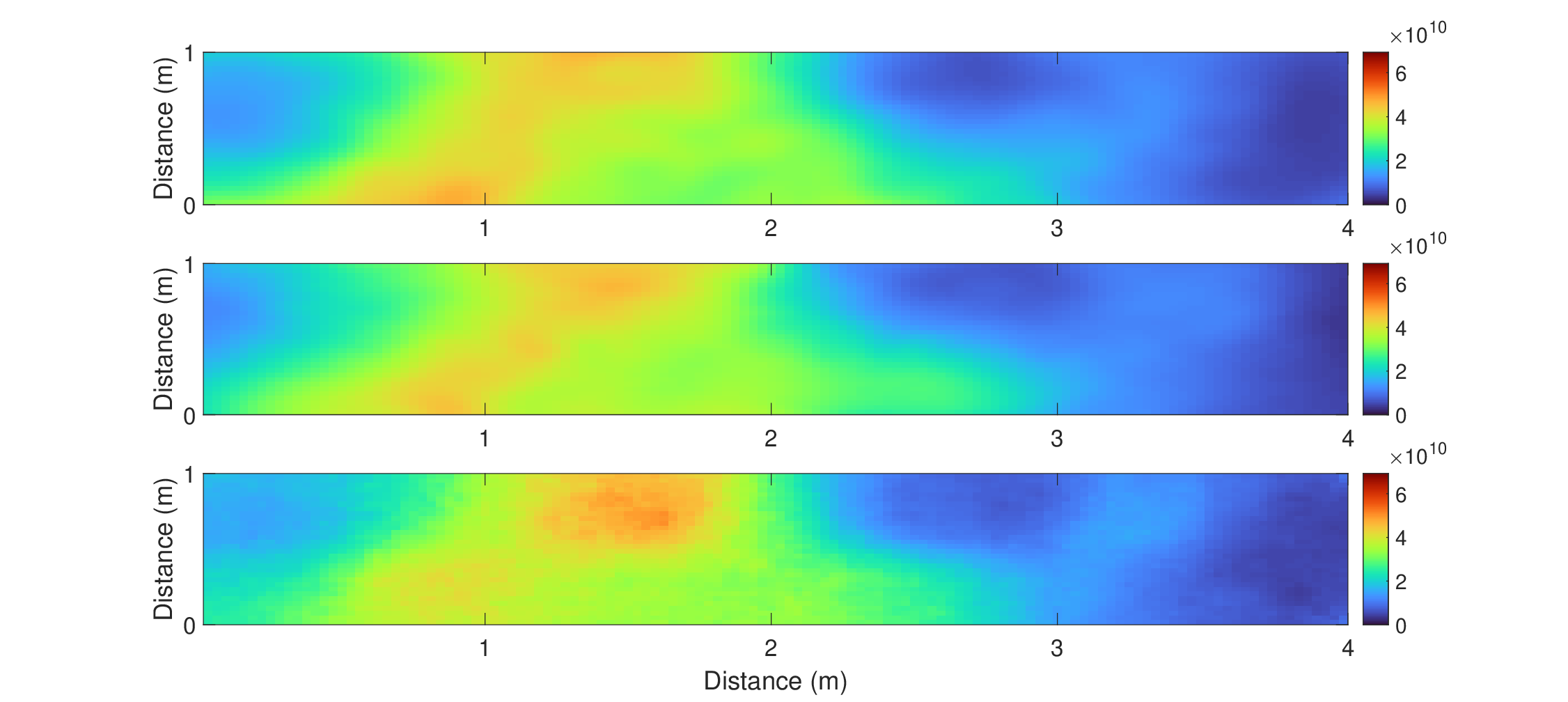}
\caption{Posterior mean values of the Young's modulus field, in Pa. Top to bottom: KL, Wavelet and LAS methods.}
\label{fig:postmean_priorsamp}
\end{figure}

Figure \ref{fig:postmean_priorsamp} shows the posterior mean field obtained with each representation.

Here we see that the three methods recover a very similar posterior mean. The small differences between the fields are likely due to the different approaches placing slightly different prior distributions on the field when they are truncated after a finite number of terms. The overall shape of this posterior mean field is in agreement with the ground truth, although the variance within the fields is smaller. This is likely due to the relatively high uncertainty on the data and thus low weight placed on the likelihood in Bayes' rule for this experiment, which was chosen primarily to provide good mixing within the chains. As a result, the posterior is noticeably influenced by the prior and the posterior mean lies closely to the prior mean, which is a uniform field with a value of $E \approx 26.1$ GPa. 

\begin{figure}[h]
\centering
\includegraphics[width=\textwidth]{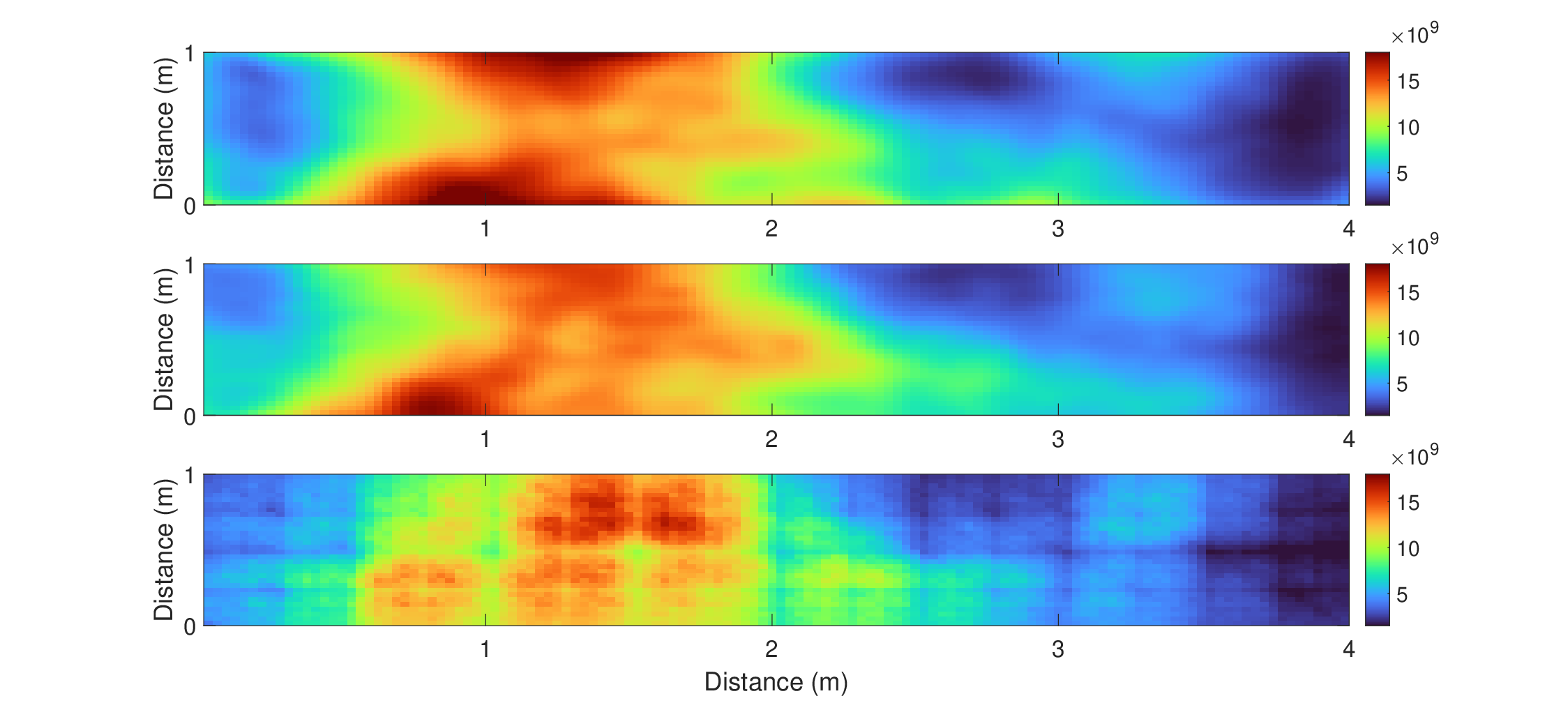}
\caption{Pointwise posterior standard deviations of the Young's modulus field, in Pa. Top to bottom: KL, Wavelet and LAS methods.}
\label{fig:postvar_priorsamp}
\end{figure}

This is also evidenced by the relatively large posterior variance. Figure \ref{fig:postvar_priorsamp} shows the pointwise posterior standard deviation of the field obtained from the finest-level samples in the multilevel algorithm. From this figure it follows that the ground truth lies within one standard deviation of the posterior mean in the majority of the domain, pointing to a rather wide but largely unbiased posterior. This figure illustrates the benefit of a full Bayesian approach: rather than just obtaining a single solution providing the best fit, we can additionally quantify the uncertainty on this estimate and specifically identify the parts of the domain where this uncertainty is greatest.

\subsubsection{Piecewise constant ground truth}

\begin{figure}[h]
\centering
\includegraphics[width=\textwidth]{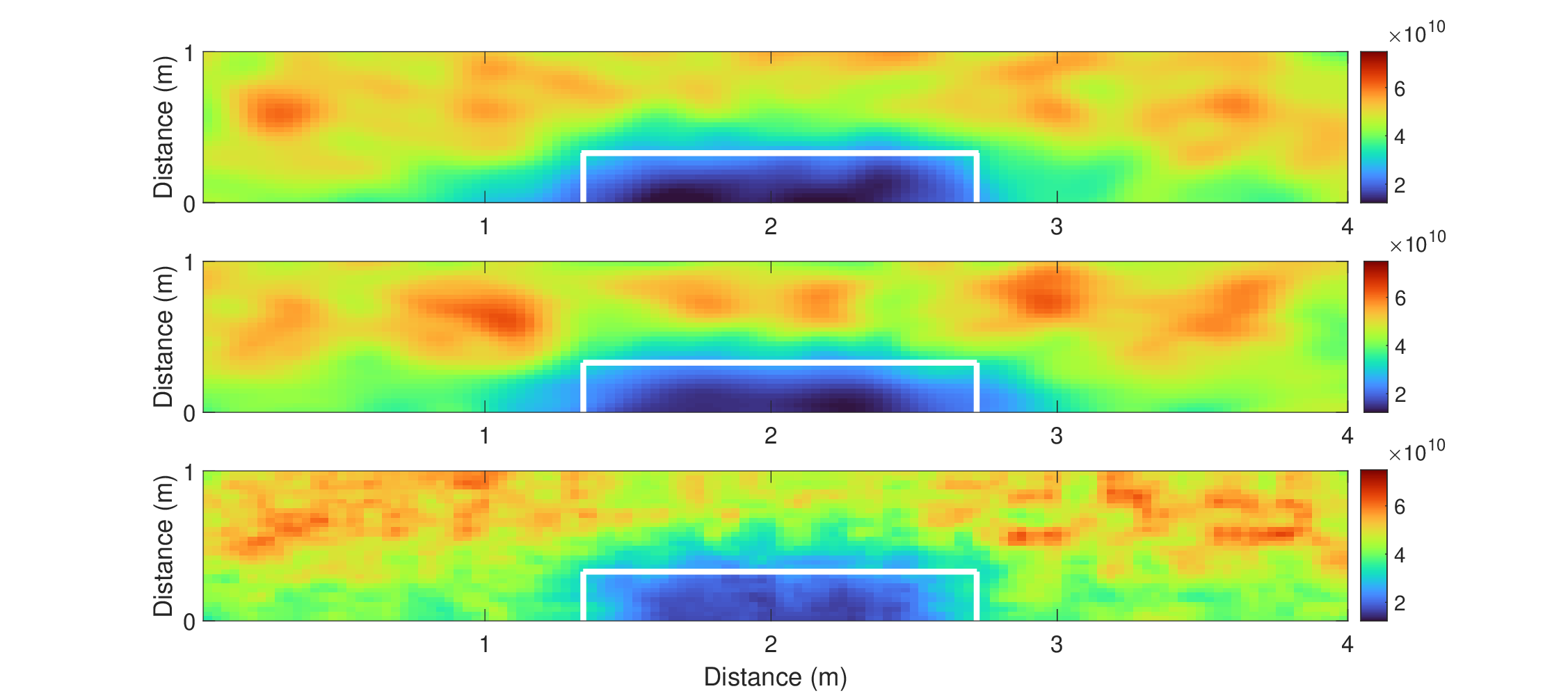}
\caption{Posterior mean values of the Young's modulus field for the second experiment, in Pa. Top to bottom: KL, Wavelet and LAS methods.}
\label{fig:postmean_piecewise}
\end{figure}

Figure \ref{fig:postmean_piecewise} shows the posterior mean Young's modulus fields for all three approaches for the piecewise constant experiment. This time, there is a considerably stronger agreement with the ground truth due to the smaller value of $\sigma_F$ in this experiment. All three methods recover the field quite well, though the posterior mean LAS field slightly misjudges the exact dimensions and value of the region of reduced $E$.

\begin{figure}[h]
\centering
\includegraphics[width=\textwidth]{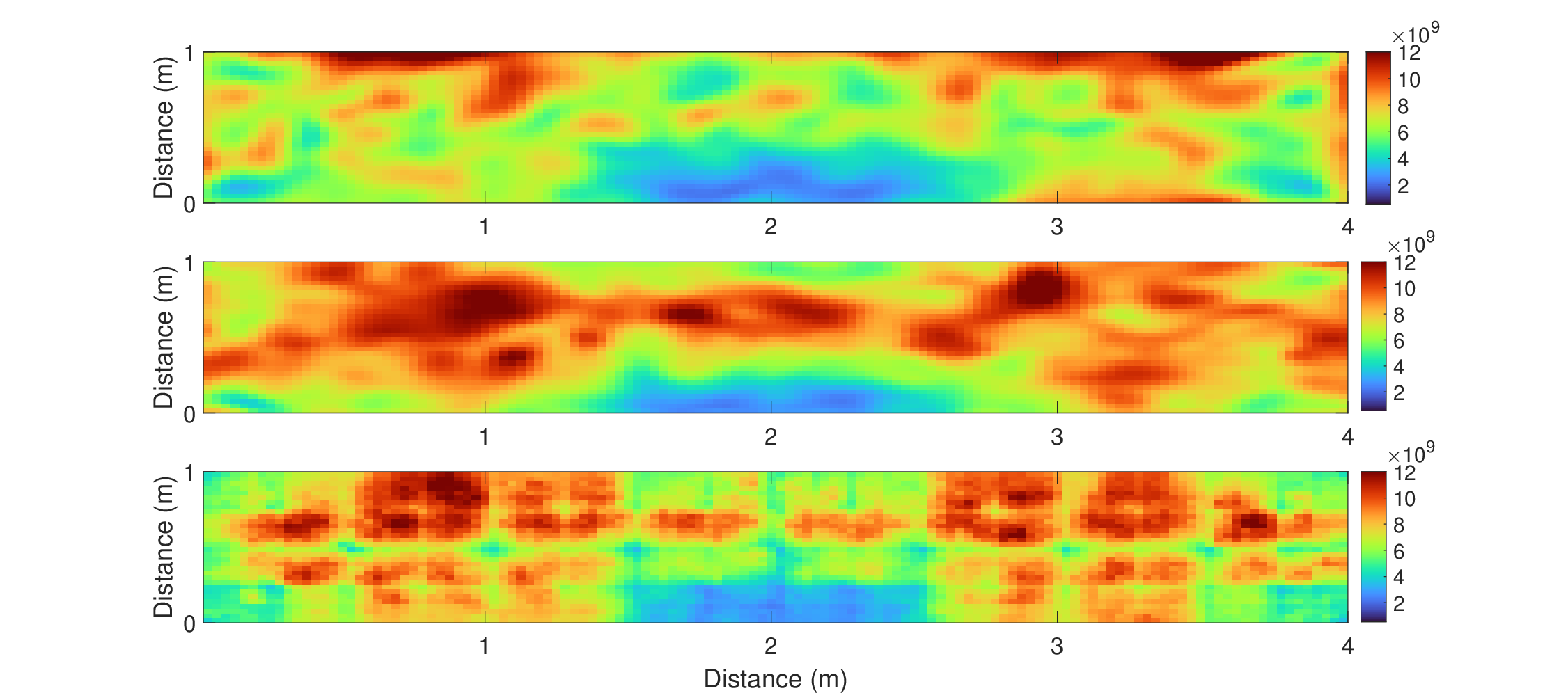}
\caption{Pointwise posterior standard deviations of the Young's modulus field for the second experiment, in Pa. Top to bottom: KL, Wavelet and LAS methods.}
\label{fig:postvar_piecewise}
\end{figure}

We can also look at the pointwise posterior variance. Figure \ref{fig:postvar_piecewise} shows the standard deviation obtained from the finest-level samples. We clearly see that all three methods more or less identify the same regions of the field where posterior uncertainty is greatest with only small differences between them. Based on these two experiments, we note that it is safe to assume that all three methods will provide users with qualitatively similar posterior estimates.

\subsection{Cost comparison}
\label{sec:cost_comparison}

\begin{figure}
\centering
\includegraphics[width=\textwidth]{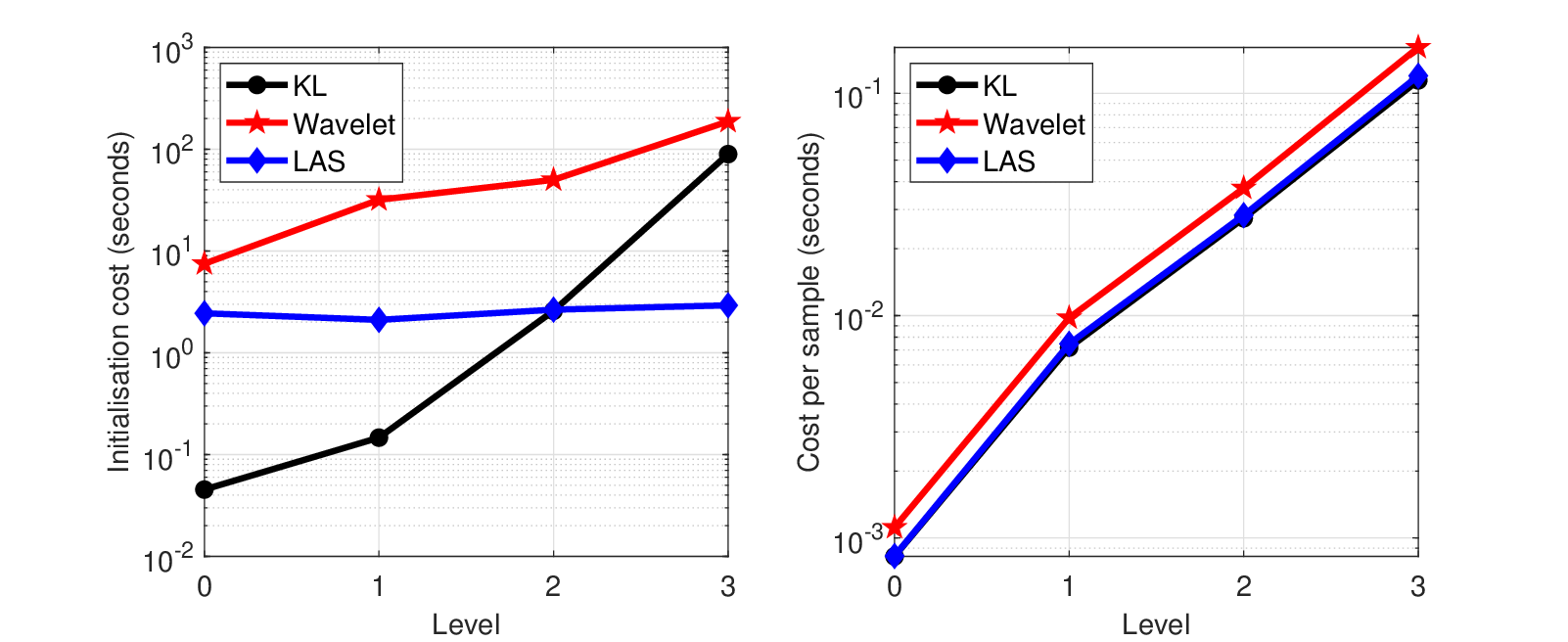}
\caption{Left: Initialisation cost to setup each of the sampling methods prior to MCMC routine. Right: Cost per sample (in seconds) at each level in the multilevel algorithm for all three representations.}
\label{fig:costs}
\end{figure}

We are additionally interested in a comparison of the three methods with respect to their computational cost. There are two aspects to the cost: initialisation cost before MCMC and sampling cost during MCMC. The initialisation step involves solving the eigenvalue problem of Equation \eqref{eq:eigenvalue_problem} for KL, generating the $b_m^{\mathbb{T}}$ in Equation \eqref{eq:waveletsum} for wavelet and calculating the coefficients $a$ and $c$ in Equations \eqref{eq:LAS_update} for LAS. The left plot of Figure \ref{fig:costs} shows this cost for each level in the multilevel hierarchy.

While the LAS method barely increases in cost for finer levels, the wavelet and especially the KL methods scale much less efficiently. It should be noted that, for the level hierarchy constructed here, even for the KL approach this initialisation cost is negligible with respect to the cost of the MCMC routine. However, if finer levels are desired, the cost of the KL method rapidly escalates due to the complexity of solving the eigenvalue problem. In these cases, this might significantly influence the total cost of the inversion routine. The LAS and to a lesser extent wavelet approaches suffer less from this.

The main cost of the algorithm is the runtime of the MCMC routine. To compare this, we will use the estimate for the sampling error $\varepsilon$ in Equation \eqref{eq:mse}. For a given computational budget, this estimate contains (i) the number of samples $N_{\ell}$ that can be generated at each level, (ii) the integrated autocorrelation time IAT$_{\ell}$ of the chain at level $\ell$ \cite{foreman-mackey_emcee_2013}, and (iii) the sample variance $\mathbb{V}(Y_{\ell})$ of the quantity of interest $Q$ at $\ell = 0$ and its corrections $Y$ at $\ell > 0$:

\begin{equation}
\varepsilon^2 = \sum_{\ell=0}^L \dfrac{\mathbb{V}(Y_{\ell})\cdot \text{IAT}_{\ell}}{N_{\ell}}.
\label{eq:mc_error}
\end{equation}

The right plot of Figure \ref{fig:costs} shows the cost per sample at each level. The KL and LAS methods have virtually the same cost per sample, while the wavelet approach is roughly 30\% more expensive per sample, mostly due to the considerably higher truncation number compared to the KL approach. 

To compare the IATs and sample variances in the different methods, we look at the two different experiments outlined above.

\subsubsection{Prior-sampled ground truth}

The right plot of Figure \ref{fig:costs} shows the estimated value of $\varepsilon$ for a given computational budget $\mathcal{C}$ for all three representations.

As the cost $\mathcal{C}$ of the multilevel MCMC algorithm scales linearly with $N_{\ell}$ it can be clearly seen from the right plot of Figure \ref{fig:costs} that all three methods exhibit the typical $ \varepsilon \simeq \mathcal{C}^{-1/2}$ convergence rate expected from the algorithm \cite{dodwell_hierarchical_2015}. However, the constants differ quite significantly, with the LAS method having a considerably lower sampling error than the other two approaches. As illustrated by Figure \ref{fig:costs}, this is not due to the individual cost per sample, as this cost is dominated by solving the FE system which occurs after constructing the realisation of $E$ and is identical for all approaches.

Instead, the differing errors are explained by looking at the factor $\mathbb{V}(Y_{\ell})\cdot \text{IAT}_{\ell}$ in Equation \eqref{eq:mc_error}. Table \ref{tab:cost_iat} lists the sample variances and IATs at the four levels of the algorithm.

\begin{table}[h]
\centering
\begin{tabular}{cc|ccc||cc|ccc}
&& \textbf{KL} & \textbf{Wavelet} & \textbf{LAS} & & & \textbf{KL} & \textbf{Wavelet} & \textbf{LAS} \\
\hline
$10^4 \cdot\mathbb{V}(Y_{\ell})$ & $\ell=0$ & 0.686 & 1.17 & 0.574 & IAT & $\ell=0$ & 1134 & 2814 & 2390 \\
& $\ell=1$ & 0.839 & 1.31 & 0.250 & & $\ell=1$ & 63 & 175 & 29 \\
& $\ell=2$ & 0.560 & 0.708 & 0.157 & & $\ell=2$ & 36 & 44 & 8 \\
& $\ell=3$ & 0.365 & 0.772 & 0.0866 & & $\ell=3$ & 9 & 15 & 4 \\
\end{tabular}
\caption{Variance of correction terms and integrated autocorrelation times of the chains at all levels, for all three field representations. Variances are multiplied by $10^4$ for ease of notation.}
\label{tab:cost_iat}
\end{table}

For the KL and wavelet models, the decay of the variances of the correction terms occurs quite slowly, with all levels having the same order of magnitude of sample variance. Only for the LAS approach is there a significant reduction of this variance for increasing levels. Regarding the autocorrelation, we additionally observe that the reduction in IAT for increasing level is comparable among all three methods, with the LAS method showing the strongest decay. This is especially notable, as the LAS approach uses considerably higher dimensional chains at all levels than the two other approaches.

\begin{figure}[h]
\centering
\includegraphics[width=\textwidth]{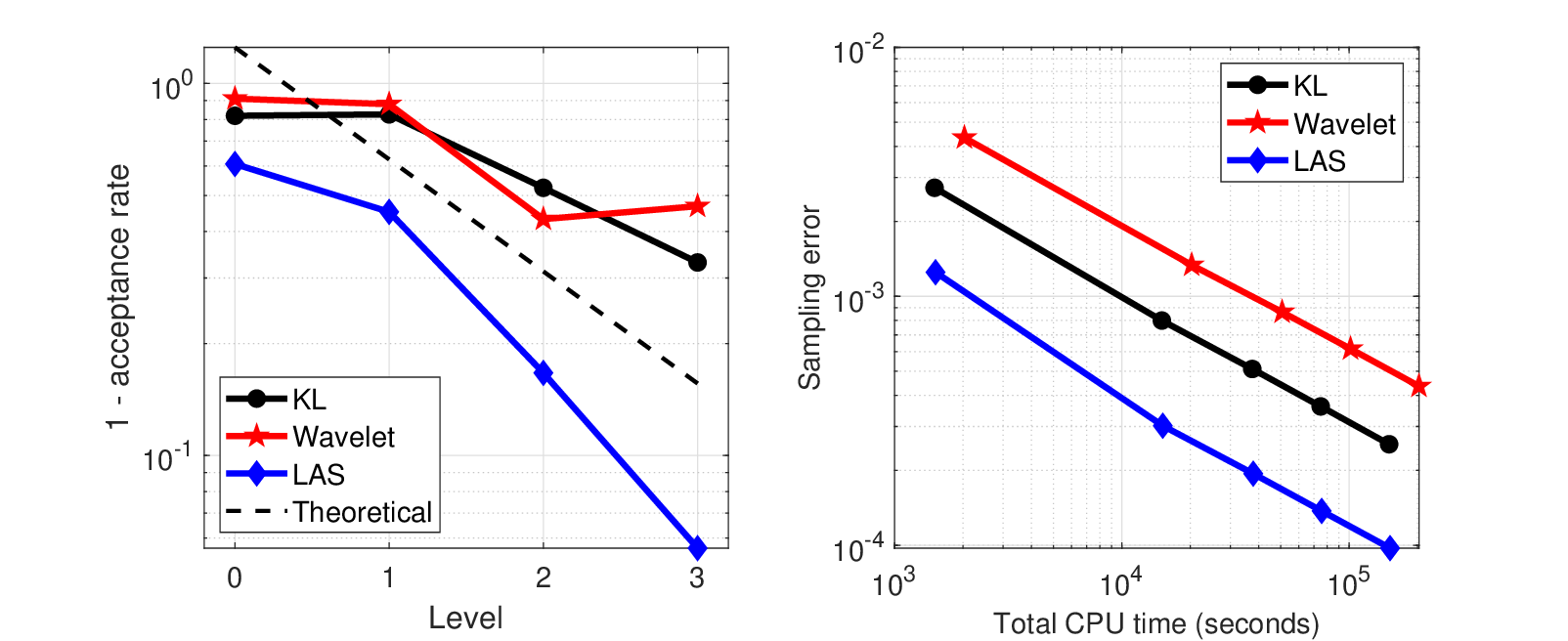}
\caption{Left: rejection rates at each level in the multilevel MCMC algorithm for all three methods for the first experiment. Right: estimates of the Monte Carlo error of Equation \eqref{eq:mc_error} for a given computational budget. Theoretical rates in dashed lines.}
\label{fig:accrates}
\end{figure}

To explain this difference in mixing, we look at the convergence of the acceptance rate of the MCMC algorithm. For acceptance rate $\alpha_{\ell}$ at each level, it is known that \cite{dodwell_hierarchical_2015}

\begin{equation}
\mathbb{E}\left[1 - \alpha_{\ell}(E_{\ell}'|E_{\ell}^n)\right] \lesssim M_{\ell-1}^a + R_{\ell-1}^{-1/d}
\end{equation}

where $R_{\ell}$ is proportional to the number of finite elements in each dimension of the $d$-dimensional grid, and $a$ is the convergence rate of the expansions from Equations \eqref{eq:kl_convrate}-\eqref{eq:las_convrate}, i.e. $a=-\nu/d$ for KL and wavelet and $a=-1/d$ for LAS. For the resolutions and truncations chosen, it holds that $R_{\ell-1}^{-1/d} \geq M_{\ell-1}^a$ for all three methods, so the convergence rate follows the discretization of the grid. The left plot of Figure \ref{fig:accrates} shows the rejection rates of all three methods, with the theoretical rate in dashed lines. The KL method more or less follows this rate.

The wavelet method does not converge as cleanly, which we believe to be likely a pre-asymptotic effect, as the convergence of the wavelet expansion is asymptotically of order $\nu/d$ but does not follow this rate smoothly as evidenced in Figure \ref{fig:truncations}.

The acceptance rate of the LAS method on the other hand converges noticeably faster than theoretically expected, which is why it results in the substantially smaller values for the IAT. We currently have no explanation for this behaviour, though we suspect it might arise from the fact that unlike the other two approaches, the LAS method ensures that the average of the subdivided elements remains the same as the value in the parent element. We believe this might cause stronger correlation of individual terms in the correction and hence a stronger decrease in sample variance of the correction terms.

\subsubsection{Piecewise constant ground truth}

We repeat the analysis for the piecewise constant ground truth experiment. First, Table \ref{tab:cost_iat2} again lists the variance of the correction terms and the IATs of the different chains. Due to the smaller step size, the IATs are far larger in this case, which will result in larger statistical errors overall. We see similar behaviour as before, with all methods having a strong reduction in IAT for levels $\ell>0$. Again the wavelet approach has the largest variance in the correction term estimate and LAS the lowest, though now the variance reduction is more pronounced in the KL approach.

\begin{table}[h]
\centering
\begin{tabular}{cc|ccc||cc|ccc}
&& \textbf{KL} & \textbf{Wavelet} & \textbf{LAS} & & & \textbf{KL} & \textbf{Wavelet} & \textbf{LAS} \\
\hline
$10^4 \cdot\mathbb{V}(Y_{\ell})$ & $\ell=0$ & 0.215 & 0.482 & 0.189 & IAT & $\ell=0$ & 4614 & 14256 & 36539 \\
& $\ell=1$ & 0.380 & 0.540 & 0.0346 & & $\ell=1$ & 209 & 291 & 372 \\
& $\ell=2$ & 0.149 & 0.272 & 0.0279 & & $\ell=2$ & 95 & 98 & 79 \\
& $\ell=3$ & 0.0858 & 0.271 & 0.0434 & & $\ell=3$ & 25 & 22 & 19 \\
\end{tabular}
\caption{Variance of correction terms and integrated autocorrelation times of the chains at all levels, for all three field representations in the second experiment. Variances are multiplied by $10^4$ for ease of notation.}
\label{tab:cost_iat2}
\end{table}

\begin{figure}[h]
\centering
\includegraphics[width=\textwidth]{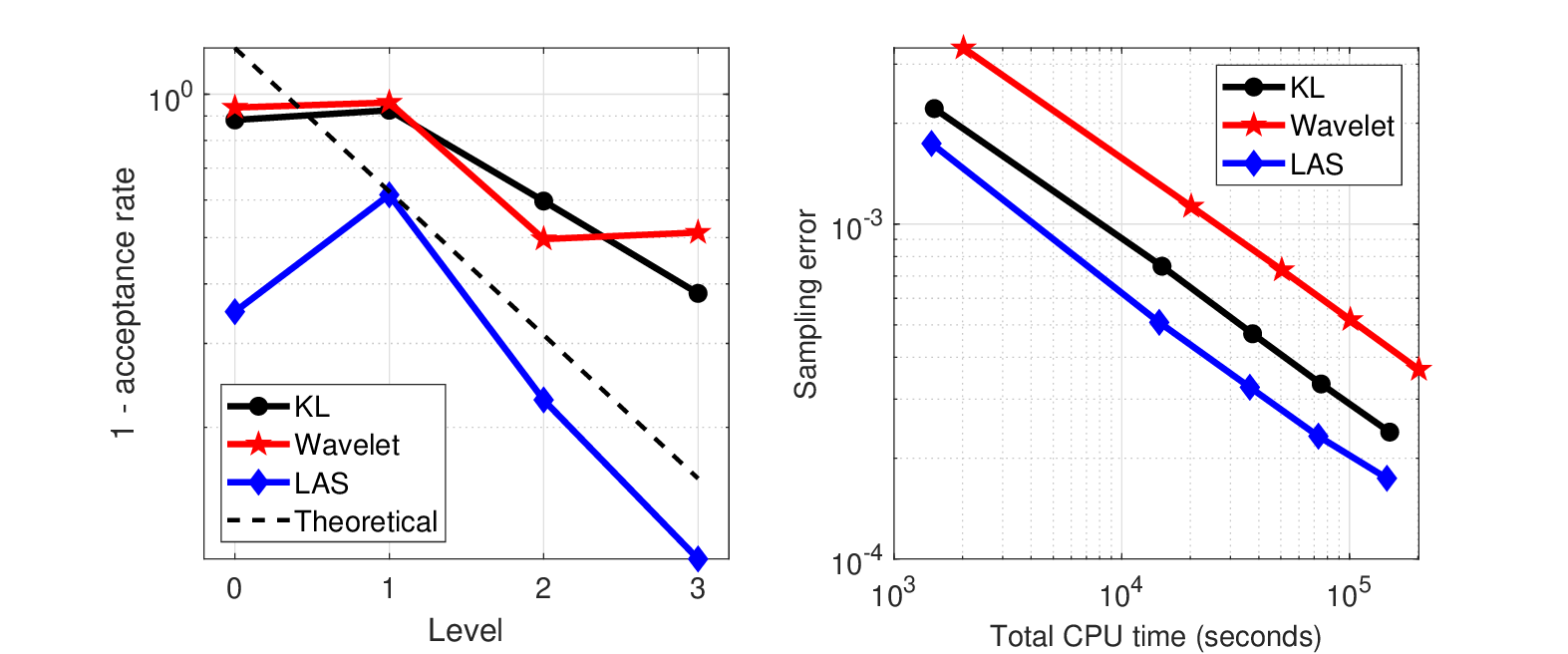}
\caption{Left: convergence of rejection rate in the multilevel MCMC algorithm for all three methods for the second experiment, with theoretical rate in dashed lines. Right: estimates of the Monte Carlo error of Equation \eqref{eq:mc_error} for a given computational budget.}
\label{fig:accrates2}
\end{figure}

Finally, Figure \ref{fig:accrates2} shows the rejection rate convergence in the left plot and the sampling error in the right plot. We again see very similar behaviour to before, with all three methods having $\mathcal{C} \simeq \varepsilon^{-2}$ rates. The rejection rates are also again quite similar, with the KL and wavelet following the theoretical convergence rate and the LAS method converging slightly faster than expected.

\section{Conclusion}
\label{ch:conclusion}

We compare three different strategies of random field sampling in terms of practical use for uncertainty quantification: the Karhunen-Lo\`eve (KL) expansion, a wavelet expansion and the method of local average subdivision (LAS). Each of these methods models random fields by constructing an expansion based on the covariance structure of a given Gaussian field, which is subsequently fed through a transformation to obtain the desired qualities of the random field. This expansion takes standard normal coefficients to model different realisations of the field, making the approach well suited for Bayesian inversion, as a standard normal prior can be used.

We use the multilevel Markov chain Monte Carlo methodology to perform a comparison between the three methods in a high-dimensional, high-resolution setting. The pCN proposal distribution is used to overcome the challenge of high-dimensionality in the chains. The three methods result in similar, only slightly different posterior estimates. Such slight differences are to be expected using three different approaches. 

Looking at the Monte Carlo error for each method, all three methods exhibit the classical $\mathcal{C} \simeq \varepsilon^{-2}$ behaviour expected from the multilevel MCMC algorithm. However, the hidden constant is found to differ quite significantly between the methods. The LAS approach consistently performs better in this regard, mainly due to its faster convergence of the acceptance rate resulting in better mixing of the Markov chains, and the smaller magnitude of the correction terms in the multilevel estimator. For these parameters, the LAS convergence appears even faster than theoretically expected. We currently have no explanation for this faster-than-expected convergence and merely note the observation.

For large-scale FE simulations with spatially heterogeneous materials, the choice of stochastic field discretization has a decisive impact on computational efficiency, even when posterior accuracy is comparable. Computationally, the LAS method exhibits better mixing than the wavelet and KL approaches at the fine levels of the multilevel algorithm. This results in a cost-to-error relation which has the same asymptotic rate for all methods, but a lower constant factor for the LAS method. The considerably larger dimensions used in the LAS expansion pose no issue due to the use of the pCN proposal distribution. Additionally, though negligible for low-resolution grids, its initialisation cost is far lower than that of the other two methods for very high-resolution grids. Its main drawback rather lies in the limited geometries for which LAS implementations currently exist. Between the other two methods, we observe that the KL method performs better than the wavelet approach, again with the same convergence rate but a smaller constant.

\section*{Acknowledgements}

The authors would like to thank Iason Papaioannou and Aretha Teckentrup for their interesting discussions on the context and interpretation of the results in this paper, as well as for their helpful comments in improving this text.

The work of the first author was supported by an SB Ph.D. fellowship of the Research Foundation Flanders (FWO) through grant 1SD1823N. Additional funding was provided by the KU Leuven Research Fund under grant C14/23/098.


\end{document}